\documentclass[final,3p,11pt,sort&compress]{elsarticle}

\usepackage{amsmath}
\usepackage{amsfonts}
\usepackage{amssymb}
\usepackage{amsthm}
\usepackage{siunitx}
\usepackage{graphicx} 
\usepackage{xspace} 
\usepackage{bm} 
\usepackage{threeparttable} 
\usepackage{subfig}
\usepackage[ruled,lined,linesnumbered]{algorithm2e}
\usepackage[english]{babel}
\usepackage{amsmath}
\usepackage{graphicx}
\usepackage[justification=centering]{caption}
\usepackage{booktabs}
\usepackage{placeins}
\usepackage{soul}
\usepackage[table]{xcolor}

\usepackage{hyperref}
\hypersetup{
 colorlinks=true,
 citecolor=blue,
 linkcolor=blue,
 urlcolor=blue}
\usepackage{tikz}
\usetikzlibrary{matrix,backgrounds,shapes}
\usetikzlibrary{positioning}
\usetikzlibrary{fit}
\usetikzlibrary{plotmarks}
\usetikzlibrary{decorations.pathreplacing}
\usepackage{pgfplots}
\usetikzlibrary{shapes}
\usetikzlibrary{positioning,arrows}
\usetikzlibrary{fadings,shapes.arrows,shadows}
\usetikzlibrary{arrows.meta}
\definecolor{darkorange}{rgb}{1.0, 0.55, 0.0}
\definecolor{Royalblue}{rgb}{0.254,0.41,0.88}
\definecolor{royalblue}{rgb}{0.254,0.41,0.88}
\usepackage{xcolor}

\newtheorem{definition}{Definition}
\newtheorem{remark}{Remark}

\usepackage{xspace}
\usepackage{amsmath}
\usepackage{amsfonts}
\usepackage{amssymb}
\usepackage{bm}
\usepackage{xspace}
\usepackage{afterpage}
\usepackage{color}
\usepackage{soul}
\usepackage{cancel}
\usepackage{gensymb}
\usepackage{tikz}
\usetikzlibrary{matrix,backgrounds,shapes}
\usetikzlibrary{positioning}
\usetikzlibrary{fit}
\usetikzlibrary{plotmarks}
\usetikzlibrary{decorations.pathreplacing}
\usepackage{pgfplots}
\usetikzlibrary{shapes}
\usetikzlibrary{positioning,arrows}
\usetikzlibrary{fadings,shapes.arrows,shadows}  
\usetikzlibrary{arrows.meta}
\definecolor{darkorange}{rgb}{1.0, 0.55, 0.0}
\definecolor{do}{rgb}{1.0, 0.55, 0.0}
\definecolor{Royalblue}{rgb}{0.254,0.41,0.88}
\usepackage{xcolor}

\definecolor{darkorange}{rgb}{1.0, 0.55, 0.0}
\definecolor{royalblue}{rgb}{0.254,0.41,0.88}

\newcommand{\Tr}{\ensuremath{^{\mr{T}}}}

\newcommand{\mr}[1]{\ensuremath{\mathrm{#1}}}
\newcommand{\fnc}[1]{\ensuremath{\mathcal{#1}}}
\newcommand{\bfnc}[1]{\ensuremath{\bm{\mathcal{#1}}}}
\newcommand{\mat}[1]{\ensuremath{\mathsf{#1}}}

\newcommand{\eg}[0]{{e.g.\@}\xspace}
\newcommand{\ie}[0]{{i.e.\@}\xspace}

\newcommand{\Th}[0]{\ensuremath{^{\mathrm{th}}}}
\newtheorem{assume}{Assumption}

\DeclareMathOperator{\diag}{diag}

\newcommand{\xm}[1]{\ensuremath{x_{#1}}}
\newcommand{\xil}[1]{\ensuremath{\xi_{#1}}}
\newcommand{\alphal}[1]{\ensuremath{\alpha_{#1}}}
\newcommand{\betal}[1]{\ensuremath{\beta_{#1}}}
\newcommand{\Nl}[1]{\ensuremath{N_{#1}}}
\newcommand{\bxil}[1]{\ensuremath{\bm{\xi}_{#1}}}
\newcommand{\bxili}[2]{\ensuremath{\bm{\xi}_{#1}^{(#2)}}}

\newcommand{\Q}[0]{\ensuremath{\bm{\fnc{Q}}}}
\newcommand{\Jk}[0]{\ensuremath{\fnc{J}}}

\newcommand{\Jdxildxm}[2]{\ensuremath{\Jk\frac{\partial\xil{#1}}{\partial\xm{#2}}}}
\newcommand{\Fxm}[1]{\ensuremath{\bm{\fnc{F}}_{\xm{#1}}}}
\newcommand{\FxmI}[1]{\ensuremath{\bm{\fnc{F}}_{\xm{#1}}^{(I)}}}
\newcommand{\FxmV}[1]{\ensuremath{\bm{\fnc{F}}_{\xm{#1}}^{(V)}}}

\newcommand{\Um}[1]{\ensuremath{\fnc{U}_{#1}}}
\newcommand{\E}[0]{\ensuremath{\fnc{E}}}

\newcommand{\GB}[0]{\ensuremath{\bm{\fnc{G}}^{(B)}}}
\newcommand{\Gzero}[0]{\ensuremath{\bm{\fnc{G}}^{(0)}}}

\newcommand{\DxiloneD}[1]{\ensuremath{\mat{D}_{\xil{#1}}^{(1D)}}}
\newcommand{\PxiloneD}[1]{\ensuremath{\mat{P}_{\xil{#1}}^{(1D)}}}
\newcommand{\QxiloneD}[1]{\ensuremath{\mat{Q}_{\xil{#1}}^{(1D)}}}
\newcommand{\SxiloneD}[1]{\ensuremath{\mat{S}_{\xil{#1}}^{(1D)}}}
\newcommand{\ExiloneD}[1]{\ensuremath{\mat{E}_{\xil{#1}}^{(1D)}}}
\newcommand{\txilalpha}[1]{\ensuremath{\bm{t}_{\alphal}}}
\newcommand{\txilbeta}[1]{\ensuremath{\bm{t}_{\betal}}}

\newcommand{\Imat}[1]{\ensuremath{\mat{I}_{#1}}}
\newcommand{\M}[0]{\ensuremath{\mat{P}}}

\newcommand{\Dxil}[1]{\ensuremath{\mat{D}_{\xil{#1}}}}
\newcommand{\Qxil}[1]{\ensuremath{\mat{Q}_{\xil{#1}}}}

\newcommand{\matFxm}[3]{\ensuremath{\mat{F}_{\xm{#1}}\left(#2,#3\right)}}

\newcommand{\qk}[1]{\ensuremath{\bm{q}_{#1}}}

\newcommand{\ones}[1]{\ensuremath{\bm{1}_{#1}}}

\newcommand{\wk}[1]{\ensuremath{\bm{w}_{#1}}}

\newcommand{\matJk}[1]{\ensuremath{{\color{orange}\mat{J}}_{#1}}}

\newcommand{\matAlmk}[3]{\ensuremath{\left[{\color{blue}\fnc{J}\frac{\partial\xil{#1}}{\partial\xm{#2}}}\right]_{#3}}}
\newcommand{\matAlmkAnal}[3]{\ensuremath{\left[{\color{red}\fnc{J}\frac{\partial\xil{#1}}{\partial\xm{#2}}}\right]_{#3}}}

\newcommand{\Ok}[0]{\ensuremath{\Omega_{\kappa}}}
\newcommand{\pOk}[0]{\ensuremath{\partial\Omega_{\kappa}}}

\newcommand{\Ohatk}[0]{\ensuremath{\hat{\Omega}_{\kappa}}}

\newcommand{\pOhatk}[0]{\ensuremath{\partial\hat{\Omega}_{\kappa}}}

\newcommand{\uk}[0]{\ensuremath{\bm{u}_{\kappa}}}
\newcommand{\ukp}[0]{\ensuremath{\bm{u}_{\kappa'}}}

\newcommand{\Cij}[2]{\ensuremath{\mat{C}_{#1,#2}}}
\newcommand{\Chatij}[2]{\ensuremath{\hat{\mat{C}}_{#1,#2}}}

\newcommand{\eonel}[1]{\ensuremath{\bm{e}_{1}}}
\newcommand{\eNl}[1]{\ensuremath{\bm{e}_{N}}}
\newcommand{\bmxi}[1]{\ensuremath{\bm{\xi}^{(#1)}}}
\newcommand{\U}[0]{\ensuremath{\fnc{U}}}

\newcommand{\MM}[0]{\ensuremath{\mat{M}}}

\newcommand{\GBlin}[0]{\ensuremath{{\fnc{G}}^{(B)}}}
\newcommand{\Gzerolin}[0]{\ensuremath{{\fnc{G}}^{(0)}}}

\newcommand{\am}[1]{\ensuremath{a_{#1}}}
\newcommand{\Bm}[1]{\ensuremath{b_{#1}}}
\newcommand{\Chatla}[2]{\ensuremath{\hat{\fnc{C}}_{#1,#2}}}
\newcommand{\Thetaa}[1]{\ensuremath{\Theta_{#1}}}

\newcommand{\matChatlaAnal}[2]{\ensuremath{\left[{\color{red}\Chatla{#1}{#2}}\right]}}

\newcommand{\thetaa}[1]{\ensuremath{\bm{\theta}_{#1}}}
\newcommand{\IP}[0]{\ensuremath{\bm{I}_{P}}}

\renewcommand{\equiv}{=}

\begin{document}

\begin{frontmatter}

\title{Optimized geometrical metrics satisfying free-stream preservation}

\author[KAUST]{Irving Reyna Nolasco\fnref{fn1}}
\ead{irvingenrique.reynanolasco@kaust.edu.sa}
\author[KAUST]{Lisandro Dalcin\fnref{fn2}} 
\ead{dalcinl@gmail.com}
\author[nia,nasa]{David C.~Del Rey Fern\'andez\fnref{fn2}}
\ead{dcdelrey@gmail.com}
\author[KAUST]{Stefano Zampini\fnref{fn2}}
\ead{stefano.zampini@kaust.edu.sa}
\author[KAUST]{Matteo Parsani\fnref{fn4}\corref{cor1}}
\ead{matteo.parsani@kaust.edu.sa}

\cortext[cor1]{Corresponding author}

\fntext[fn1]{Ph.D. student}
\fntext[fn2]{Research Scientist}
\fntext[fn4]{Assistant Professor}

\address[KAUST]{King Abdullah University of Science and Technology (KAUST), 
  Computer Electrical and Mathematical Science and Engineering Division (CEMSE), 
  Extreme Computing Research Center (ECRC), Thuwal, Saudi Arabia}
\address[nia]{National Institute of Aerospace, Hampton, Virginia, United States}
\address[nasa]{Computational AeroSciences Branch, NASA Langley Research Center,
  Hampton, Virginia, United States}

\begin{abstract}
Computational fluid dynamics and aerodynamics, which complement more expensive 
empirical approaches, are critical for developing aerospace vehicles. During 
the past three decades, computational aerodynamics capability has improved 
remarkably, following advances in computer hardware and algorithm development.  
However, for complex applications, the demands on computational fluid dynamics continue to increase
 in a quest to gain a few percent improvements in 
accuracy. Herein, we numerically demonstrate that optimizing the metric terms
which arise from smoothly mapping each cell to a reference element,
lead to a solution whose accuracy is 
practically never worse and often noticeably better than the one obtained using 
the widely adopted Thomas and Lombard metric terms computation 
(Geometric conservation law and its application to flow computations on moving grids, AIAA Journal, 1979).
Low and high-order accurate entropy stable
schemes on distorted, high-order tensor product elements are used to 
simulate three-dimensional inviscid and viscous compressible test cases for which an analytical
solution is known.
\end{abstract}


\end{frontmatter}

\section{Introduction}
In recent years, with the continuous growth
of computing capability 
and in an effort to achieve more accurate numerical simulations 
over a broad class of engineering problems,
computational fluid dynamics (CFD) has gradually shifted 
towards high-order accurate simulations; 
see, for instance, \cite{wang_hom_review_2012,huynh_hom_review_2014,abgrall_hom_hal_2017}.
Modern unstructured high-order methods 
include discontinuous Galerkin (DG), 
spectral difference (SD), and flux reconstruction (FR) methods and can produce 
highly accurate solutions with 
minimum numerical dispersion and dissipation. Although DG, SD, and FR methods are well suited for 
smooth solutions, numerical instabilities may occur if the flow contains 
discontinuities or under-resolved physical features. A variety of mathematical 
stabilization strategies are commonly used to alleviate these issues (e.g., 
filtering \cite{hesthaven_nodal_dg_2008}, weighted essentially non-oscillatory schemes \cite{zhu_dg_weno_2013}, 
artificial dissipation, over-integration, and limiters); however, the use of such techniques for practical complex flow applications often times rely on heuristic (\eg, tunable parameters) or 
results in schemes that are not always stable (e.g., over-integration)

For the compressible Navier--Stokes equations, a very promising and mathematically rigorous alternative consists in constructing
discrete operators that are non-linearly stable (\ie, entropy stable) and simultaneously conserve 
mass, momentum, and energy, while enforcing a secondary entropy constraint. This 
strategy consists in first identifying a non-linear neutrally stable flux for the compressible Euler
equations, and then adding an appropriate amount of dissipation in order to achieve entropy stability at the semi-discrete level, 
 thus enhancing the stability of the base operator
\cite{fisher_entropy_stability_fd_2013,parsani_entropy_stability_solid_wall_2015,pazner_es_line_dg_2019,comparison}. 
In this work, we build our study on
conforming entropy stable discontinuous collocation methods, 
constructed by using the summation-by-parts (SBP) operators and the simultaneous-approximation-terms
(SAT) framework 
\cite{carpenter_ssdc_2014,parsani_entropy_stability_solid_wall_2015,carpenter_entropy_stability_ssdc_2016,parsani_ssdc_staggered_2016}. However, 
the proposed methodology can be immediately applied to any of the aforementioned
spatial discretization techniques.

In CFD, simulations in complex geometries are performed
on the union of piece-wise smooth sub-domains, also called elements or cells,
that are smoothly mapped to a reference element. 
In this element, the derivative terms appearing in the system
of partial differential equations (PDEs) are actually evaluated.
In the mapped system of 
equations, the fluxes include the Jacobians of the transformation and the transformation metrics, which depend on derivatives 
of the transformation. On the one hand, when the element sides are straight, the mapping is linear in each coordinate direction and the
metrics are constant. On the other hand, when the element boundaries are curved, these metric terms vary across the element. 

At the continuous level,
the metric terms naturally satisfy a set of identities \cite{vinokur_euler_curvilinear_coordinates_1974,thomas_gcl_1979,thompson_boundary_fitted_coordinate_1982}, and the importance of satisfying these identities at the discrete level
has long been recognized; see, for instance, 
\cite{thompson_boundary_fitted_coordinate_1982,hirsch_cfd_volume_2_1990,vinokur_fd_fv_1989}.
One of the consequences 
is that a constant 
free-stream solution is exactly preserved for all time, independent of the chosen coordinate system. 
Failure to preserve the free-stream condition frequently leads 
to spurious source terms that introduce errors in the solution and 
can be catastrophic in many applications.

The idea of approximating the metric terms
in such a way that certain physical quantities are preserved goes back to the early days of 
CFD. The terminology ``geometric conservation law'' (GCL) 
was introduced in 1979 by Thomas and Lombard \cite{thomas_gcl_1979}, who found 
that finite difference approximations that satisfied the metric 
identities in two dimensions failed in three dimensions. Such observation led 
to the re-write of the metric terms in a ``conservative form'', which,
when approximated with central differences, satisfy the metric identities,
and thus lead to free-stream preservation. This concept was subsequently 
extended to geometrically characterize conservative numerical schemes as algorithms that
preserve the entire state of a uniform flow.

In previous works \cite{fernandez_entropy_stable_p_ref_nasa_2019,fernandez_entropy_stable_p_euler_2019,fernandez_entropy_stable_hp_ref_snpdea_2019},
the metric terms are constructed for non-conforming discretization at the cell interfaces by
solving a strictly convex quadratic optimization problem \cite{crean_entropy_stable_sbp_curvilinear_euler}, whose 
solution guarantees
entropy conservation and 
free-stream preservation.
Herein, we numerically show that optimizing the metric terms as proposed in \cite{fernandez_entropy_stable_p_ref_nasa_2019,fernandez_entropy_stable_p_euler_2019,fernandez_entropy_stable_hp_ref_snpdea_2019}
leads, even for conforming interfaces, to a solution whose accuracy is 
practically never worse and often noticeably better than the one obtained using the widely adopted Thomas and Lombard approach \cite{thomas_gcl_1979}.
Thus, we conclude that the pre-processing step of optimizing the metric terms 
 can be used in a computational framework as a unique and viable approach for conforming and $h/p$ non-conforming
 interfaces. In addition, this choice greatly simplifies the solver and allows important code re-utilization.

The paper is organized as follows. In Section \ref{sec:notation} we introduce the
notation used in this work. The coordinate transformation from physical
to computational space and the key constraints that have to be satisfied by the discrete
metric terms are introduced in the context of the linear convection-diffusion equation
in Section \ref{sec:adv_diff_eq}. In the same section, the metric solution mechanics
for conforming interior and boundary faces are also presented. Section \ref{sec:compressible_nse} deals with 
the compressible Navier--Stokes equations 
and its semi-discretization
using entropy stable SBP-SAT operators of any order.
Numerical results for three-dimensional inviscid and viscous test cases for which an analytical solution is known are presented in Section \ref{sec:numerical_results}.
Simulations are performed using low and high-order accurate entropy stable
schemes on distorted, high-order tensor product elements \cite{carpenter_ssdc_2014,parsani_entropy_stability_solid_wall_2015,carpenter_entropy_stability_ssdc_2016}.
Conclusions are drawn in Section \ref{sec:conclusions}.

\section{Notation}\label{sec:notation}
The notation used in this work has been adopted from \cite{fernandez_entropy_stable_p_euler_2019}. 
Partial differential equations
(PDEs) are discretized on tensor-product cells having Cartesian computational coordinates denoted by 
the triple $(\xil{1},\xil{2},\xil{3})$, where the physical coordinates are denoted by the triple 
$(\xm{1},\xm{2},\xm{3})$. Vectors are represented by lowercase bold font, for example $\bm{u}$, 
while matrices are represented using sans-serif font, for example, $\mat{B}$. Continuous 
functions on a space-time domain are denoted by capital letters in script font.  For example, 
\begin{equation*}
\fnc{U}\left(\xil{1},\xil{2},\xil{3},t\right)\in L^{2}\left(\left[\alphal{1},\betal{1}\right]\times
\left[\alphal{2},\betal{2}\right]\times\left[\alphal{3},\betal{3}\right]\times\left[0,T\right]\right)
\end{equation*}
represents a square integrable function, where $t$ is the temporal coordinate. The restriction of such 
function onto a set of mesh nodes is denoted by lower case bold font. For example, the restriction of 
$\fnc{U}$ onto a grid of $\Nl{1}\times\Nl{2}\times\Nl{3}$ nodes is given by the vector
\begin{equation*}
\bm{u} = \left[\fnc{U}\left(\bxili{}{1},t\right),\dots,\fnc{U}\left(\bxili{}{N},t\right)\right]\Tr,
\end{equation*}
where $N$ is the total number of nodes ($N\equiv\Nl{1}\Nl{2}\Nl{3}$), and the square brackets are used 
to delineate vectors and matrices, as well as ranges for variables (the context will make clear which meaning is being used). Moreover, $\bm{\xi}$ is a vector of vectors 
constructed from the three vectors $\bxil{1}$, $\bxil{2}$, and $\bxil{3}$, which are 
vectors of size $\Nl{1}$, $\Nl{2}$, and $\Nl{3}$ and contain the coordinates of the mesh in 
the three computational directions, respectively. Finally, $\bxil{}$ is constructed as 
\begin{equation*}
\bxil{}(3(i-1)+1:3i)\equiv  \bxili{}{i}
\equiv\left[\bxil{1}(i),\bxil{2}(i),\bxil{3}(i)\right]\Tr,
\end{equation*}
where the notation $\bm{u}(i)$ means the $i\Th$ entry of the vector $\bm{u}$ and $\bm{u}(i:j)$ is the subvector 
constructed from $\bm{u}$ using the $i\Th$ through $j\Th$ entries (\ie, Matlab notation is used).

 Oftentimes, monomials are discussed and the following notation is used:
\begin{equation*}
\bxil{l}^{j} \equiv \left[\left(\bxil{l}(1)\right)^{j},\dots,\left(\bxil{l}(\Nl{l})\right)^{j}\right]\Tr,
\end{equation*}
with the  convention that $\bxil{l}^{j}=\bm{0}$ for $j<0$.

Herein, one-dimensional SBP operators are used to discretize derivatives. 
The definition of a one-dimensional SBP operator in the $\xil{l}$ direction, $l=1,2,3$, 
is~\cite{DCDRF2014,Fernandez2014,Svard2014}
\begin{definition}\label{SBP}
\textbf{Summation-by-parts operator for the first derivative}: A matrix operator, 
$\DxiloneD{l}\in\mathbb{R}^{\Nl{l}\times\Nl{l}}$, is an SBP operator of degree $p$ approximating the derivative 
$\frac{\partial}{\partial \xil{l}}$ on the domain $\xil{l}\in\left[\alphal{l},\betal{l}\right]$ with nodal 
distribution $\bxil{l}$ having $\Nl{l}$ nodes, if 
\begin{enumerate}
\item $\DxiloneD{l}\bxil{l}^{j}=j\bxil{l}^{j-1}$, $j=0,1,\dots,p$;
\item $\DxiloneD{l}\equiv\left(\PxiloneD{l}\right)^{-1}\QxiloneD{l}$, where the norm matrix, 
$\PxiloneD{l}$, is symmetric positive definite;
\item $\QxiloneD{l}\equiv\left(\SxiloneD{l}+\frac{1}{2}\ExiloneD{l}\right)$, 
$\SxiloneD{l}=-\left(\SxiloneD{l}\right)\Tr$, $\ExiloneD{l}=\left(\ExiloneD{l}\right)\Tr$, \\ 
$\ExiloneD{l} = \diag\left(-1,0,\dots,0,1\right)=\eNl{l}\eNl{l}\Tr-\eonel{l}\eonel{l}\Tr$, 
$\eonel{l}\equiv\left[1,0,\dots,0\right]\Tr$, and $\eNl{l}\equiv\left[0,0,\dots,1\right]\Tr$. 
\end{enumerate}
Thus, a degree $p$ SBP operator is 
one that differentiates exactly monomials up to degree $p$.
\end{definition}

In this work, one-dimensional SBP operators are extended to multiple dimensions 
using tensor products ($\otimes$).  The tensor product between the matrices $\mat{A}$ and $\mat{B}$ 
is given as $\mat{A}\otimes\mat{B}$. When referencing individual entries in a matrix the notation $\mat{A}(i,j)$ 
is used, which means 
the $i\Th$ $j\Th$ entry in the matrix $\mat{A}$.

The focus of this paper is exclusively on diagonal-norm SBP operators. Moreover, the same 
one-dimensional SBP operator is used in each direction, each operating on $N_l$ nodes. 
Specifically, diagonal-norm SBP operators constructed on the Legendre--Gauss--Lobatto (LGL) 
nodes are used, \ie, a discontinuous Galerkin collocated spectral element approach is utilized
(see, for instance, 
\cite{carpenter_ssdc_2014,parsani_entropy_stability_solid_wall_2015,carpenter_entropy_stability_ssdc_2016,parsani_ssdc_staggered_2016,gassner_entropy_shallow_water_2016,Gassner2016}).

%

When solving PDEs numerically, the physical domain $\Omega\subset\mathbb{R}^{3}$, 
with boundary $\Gamma\equiv\partial\Omega$, with Cartesian coordinates $\left(\xm{1},\xm{2},\xm{3}\right)\subset\mathbb{R}^{3}$, 
is partitioned into $K$ non-overlapping elements. The domain of the $\kappa^{\text{th}}$ element is denoted by 
$\Ok$ and has boundary $\pOk$. Numerically, we solve PDEs in computational 
coordinates, $\left(\xil{1},\xil{2},\xil{3}\right)\subset\mathbb{R}^{3}$, 
where each $\Ok$ is locally transformed to the reference element $\Ohatk$, with boundary $\pOhatk$,
using a pull-back curvilinear coordinate transformation which satisfies the following assumption:
\begin{assume}\label{assume:curv}
Each element in physical space is transformed using 
a local and invertible curvilinear coordinate transformation that is compatible at 
shared interfaces, meaning that the push-forward element-wise mappings are continuous across physical element interfaces.
Note that this is the standard assumption requiring that the curvilinear coordinate transformation is water-tight.
\end{assume}

\begin{figure}
 \begin{center}
   \includegraphics[width=0.6\textwidth]{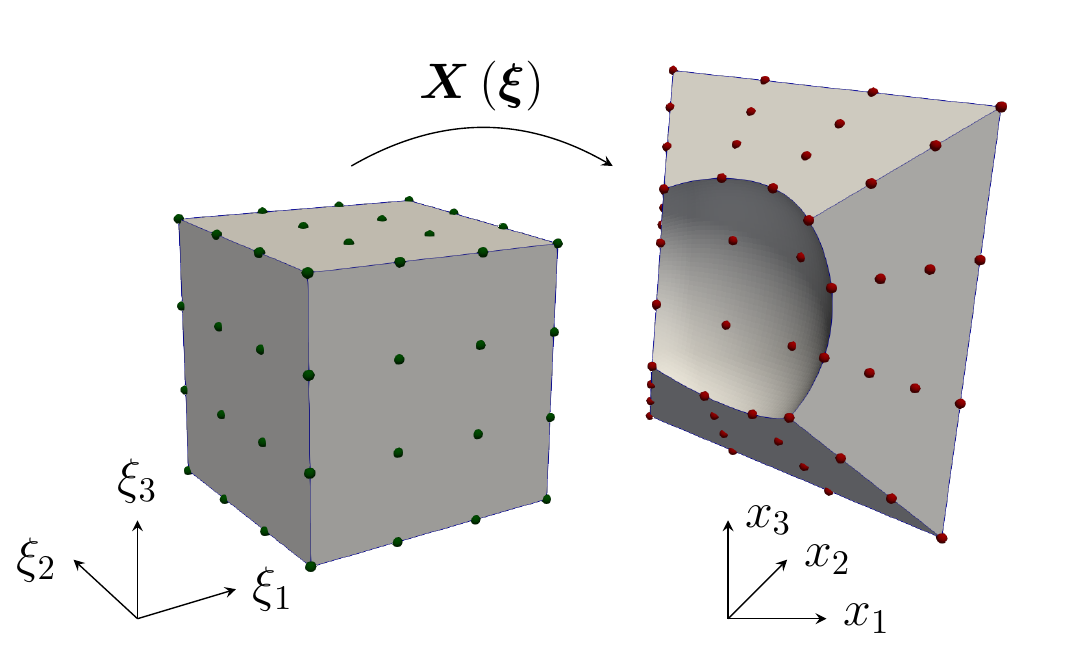}
   \caption{The reference element and its mapping to an element in the mesh.}
   \label{fig:ref_cell}
  \end{center}
\end{figure}

Precisely, one maps from the reference coordinates $\left(\xil{1},\xil{2},\xil{3}\right) \in [-1,1]^{3}$
to the physical element (see Figure \ref{fig:ref_cell}) by the push-forward transformation 
\begin{equation}\label{eq:3d_mapping}
 \left(\xm{1},\xm{2},\xm{3}\right) = X \left(\xil{1},\xil{2},\xil{3}\right),
\end{equation}
which, in the presence of curved elements, is usually a high-order degree polynomial.

\section{The linear convection-diffusion equation}\label{sec:adv_diff_eq}
Many of the technical details for constructing conservative and stable discretizations 
for the compressible Navier--Stokes, including the metric terms, are present for the discretization of the 
linear convection-diffusion equation. 
The linear convection-diffusion equation 
in Cartesian physical coordinates is given as  
\begin{equation}\label{eq:cartconvectiondiffusion}
  \begin{aligned}
    &\frac{\partial\U}{\partial t}+\sum\limits_{m=1}^{3}\frac{\partial\left(\am{m}\U\right)}{\partial\xm{m}}=
  \sum\limits_{m=1}^{3}\frac{\partial^{2}(\Bm{m}\U)}{\partial\xm{m}^{2}},
 &&
\forall \left(\xm{1},\xm{2},\xm{3}\right)\in\Omega,\quad t\ge 0,\\
  &\U\left(\xm{1},\xm{2},\xm{3},t\right)=\GBlin\left(\xm{1},\xm{2},\xm{3},t\right), && \forall \left(\xm{1},\xm{2},\xm{3}\right)\in\Gamma,\quad t\ge 0,\\
&\U\left(\xm{1},\xm{2},\xm{3},0\right)=\Gzerolin\left(\xm{1},\xm{2},\xm{3}\right), &&
  \forall \left(\xm{1},\xm{2},\xm{3}\right)\in\Omega,
  \end{aligned}
\end{equation}
where $(\am{m}\U)$ are the inviscid fluxes, $\am{m}$ are the (constant) components of the convection speed, 
$\frac{\partial(\Bm{m}\U)}{\partial\xm{m}}$ are the viscous fluxes, and $\Bm{m}$ are the (constant and positive) diffusion coefficients.  The boundary data,
$\GBlin$, and the initial condition, $\Gzerolin$, are assumed to be in $L^{2}(\Omega)$,
with the further assumption that $\GBlin$ is prescribed so that either energy conservation or energy stability is achieved.

Since derivatives are approximated with differentiation operators defined in computational space, we use the Jacobian of the push-forward mapping and the chain rule
\begin{equation*}
  \frac{\partial}{\partial\xm{m}}=
  \sum\limits_{l=1}^{3}\frac{\partial\xil{l}}{\partial\xm{m}}\frac{\partial}{\partial\xil{l}},\quad
  \frac{\partial^{2}}{\partial\xm{m}^{2}}=
  \sum\limits_{l,a=1}^{3}\frac{\partial\xil{l}}{\partial\xm{m}}
  \frac{\partial}{\partial\xil{l}}\left(
  \frac{\partial\xil{a}}{\partial \xm{m}}\frac{\partial}{\partial\xil{a}}  
  \right),
\end{equation*} 
to transform Equation \eqref{eq:cartconvectiondiffusion} from physical to computational space as
 \begin{equation}\label{eq:convectiondiffusionchain}
  \Jk\frac{\partial\U}{\partial t}
 +\sum\limits_{l,m=1}^{3}\Jk\frac{\partial\xil{l}}{\partial\xm{m}}
  \frac{\partial \left(a_{m}\U\right)}{\partial\xil{l}}=\sum\limits_{l,a,m=1}^{3}
  \Jk\frac{\partial\xil{l}}{\partial\xm{m}}\frac{\partial}{\partial\xil{l}}
  \left(\frac{\partial\xil{a}}{\partial \xm{m}}
  \frac{\partial(\Bm{m}\U)}{\partial\xil{a}}
  \right),
\end{equation}
where $\Jk$ is the determinant of the metric Jacobian.
Bringing the metric terms $\Jdxildxm{l}{m}$ inside the derivative, and using the product rule, gives
 \begin{equation}\label{eq:convectiondiffusionstrong1}
  \begin{split}
  \Jk\frac{\partial\U}{\partial t}+\sum\limits_{l,m=1}^{3}
  \frac{\partial}{\partial\xil{l}}\left(\Jdxildxm{l}{m}a_{m}\U\right)
-&\sum\limits_{l,m=1}^{3}a_{m}\U\frac{\partial}{\partial\xil{l}}\left(\Jdxildxm{l}{m}\right)
 =\\
 \sum\limits_{l,a,m=1}^{3}
  \frac{\partial}{\partial\xil{l}}
  \left(\Jk\frac{\partial\xil{l}}{\partial\xm{m}}\frac{\partial\xil{a}}{\partial \xm{m}}
  \frac{\partial(\Bm{m}\U)}{\partial\xil{a}}\right)
-&\sum\limits_{l,a,m=1}^{3}
\frac{\partial\xil{a}}{\partial \xm{m}}
  \frac{\partial(\Bm{m}\U)}{\partial\xil{a}}
  \frac{\partial}{\partial\xil{l}}\left( \Jk\frac{\partial\xil{l}}{\partial\xm{m}}\right).
  \end{split}
\end{equation}
The last terms on the left- and right-hand sides of~\eqref{eq:convectiondiffusionstrong1} are zero via the GCL relations
\begin{equation}\label{eq:GCL}
    \sum\limits_{l=1}^{3}\frac{\partial}{\partial\xil{l}}\left(\Jdxildxm{l}{m}\right)=0,\quad m=1,2,3,
\end{equation}
leading to the strong conservation form of the convection-diffusion equation in computational space
 \begin{equation}\label{eq:convectiondiffusionstrong}
  \Jk\frac{\partial\U}{\partial t}+\sum\limits_{l,m=1}^{3}
  \frac{\partial}{\partial\xil{l}}\left(\Jdxildxm{l}{m}a_{m}\U\right)= \sum\limits_{l,a,m=1}^{3}
  \frac{\partial}{\partial\xil{l}}
  \left(\Jk\frac{\partial\xil{l}}{\partial\xm{m}}\frac{\partial\xil{a}}{\partial \xm{m}}
  \frac{\partial(\Bm{m}\U)}{\partial\xil{a}}\right).
\end{equation}

Now, consider discretizing Equation~\eqref{eq:convectiondiffusionstrong} by using 
the following differentiation matrices
\begin{equation*}
\Dxil{1}\equiv\DxiloneD{1}\otimes\Imat{N_2}\otimes\Imat{N_3},\;
\Dxil{2}\equiv\Imat{N_1}\otimes\DxiloneD{2}\otimes\Imat{N_3},\;
\Dxil{3}\equiv\Imat{N_1}\otimes\Imat{N_2}\otimes\DxiloneD{3},
\end{equation*} 
where $\Imat{N_l}$ is an $N_l\times N_l$ identity matrix and $N_l$ is the number of LGL points per direction in a given element. The diagonal 
matrix containing the metric Jacobian is defined as
\begin{equation*}
  \matJk{\kappa}\equiv\diag\left(\Jk(\bmxi{1}),\dots,\Jk(\bmxi{\Nl{\kappa}})\right),
\end{equation*}
while the diagonal matrix of the metric terms, $\matAlmk{l}{m}{\kappa}$, has to be chosen to be a discretization of
\begin{equation*}
  \diag\left(\Jdxildxm{l}{m}(\bmxi{1}),\dots,
  \Jdxildxm{l}{m}(\bmxi{\Nl{\kappa}})\right),
\end{equation*}
where $\Nl{\kappa}\equiv N_1 N_2 N_3$ is the total number of nodes in element $\kappa$.
Using this nomenclature, the discretization of~\eqref{eq:convectiondiffusionstrong} on the $\kappa\Th$ element reads
 \begin{equation}\label{eq:convectionstrongdisc}
  \begin{split}
    &\matJk{\kappa}\frac{\mr{d}\uk}{\mr{d}t}+\sum\limits_{l,m=1}^{3}a_{m}\Dxil{l}\matAlmk{l}{m}{\kappa}\uk=
    \sum\limits_{l,m,a=1}^{3}\Bm{m}\Dxil{l}\matJk{\kappa}^{-1}\matAlmk{l}{m}{\kappa}\matAlmk{a}{m}{\kappa}\Dxil{a}\uk
    +\bm{SAT}_{\kappa}, 
  \end{split}
\end{equation}
where $\bm{SAT}_{\kappa}$ is the vectors of the SATs used to 
impose boundary conditions and inter-element connectivity \cite{parsani_entropy_stable_interfaces_2015,carpenter_entropy_stable_staggered_2015}. The $\bm{SAT}_{\kappa}$ vector is in general composed from inviscid and viscous contributions, \ie $\bm{SAT}_{\kappa} = \bm{SAT}^{(I)}_{\kappa} + \bm{SAT}^{(V)}_{\kappa}$.

Unfortunately, the scheme~\eqref{eq:convectionstrongdisc} is not guaranteed to be stable. However, a well-known remedy is to canonically split the inviscid terms into
one half of the inviscid terms in~\eqref{eq:convectiondiffusionchain} 
and one half of the inviscid terms in~\eqref{eq:convectiondiffusionstrong1} 
(see, for instance, \cite{carpenter_entropy_stable_staggered_2015}), while the viscous terms are treated 
in strong conservation form.  In the continuum, this process leads to
\begin{equation}\label{eq:convectiondiffusionsplit}
  \begin{split}
  &\Jk\frac{\partial\U}{\partial t}+\frac{1}{2}\sum\limits_{l,m=1}^{3}\left\{
    \frac{\partial}{\partial\xil{l}}\left(\Jdxildxm{l}{m}a_{m}\U\right)+
     \Jdxildxm{l}{m}\frac{\partial}{\partial\xil{l}}\left(a_{m}\U\right)
    \right\}\\&-\frac{1}{2}\sum\limits_{l,m=1}^{3}\left\{
    a_{m}\U\frac{\partial}{\partial\xil{l}}\left(\Jdxildxm{l}{m}\right)\right\}=
\sum\limits_{l,a,m=1}^{3}
  \frac{\partial}{\partial\xil{l}}
  \left(\Jk\frac{\partial\xil{l}}{\partial\xm{m}}\frac{\partial\xil{a}}{\partial \xm{m}}
    \frac{\partial(\Bm{m}\U)}{\partial\xil{a}}\right),  
  \end{split}
\end{equation}
where the last set of terms on the left-hand side are zero by the GCL conditions~\eqref{eq:GCL}. Then, a 
stable semi-discrete form is constructed in the same manner as the split form~\eqref{eq:convectiondiffusionsplit} by 
discretizing the inviscid portion of~\eqref{eq:convectiondiffusionchain} and~\eqref{eq:convectiondiffusionstrong} using $\Dxil{l}$, $\matJk{\kappa}$, and  
$\matAlmk{l}{m}{\kappa}$, and by averaging the results. The viscous terms result from the discretization of the viscous portion of~\eqref{eq:convectiondiffusionstrong}. 
This procedure yields
\begin{equation}\label{eq:convectionsplitdisc}
  \begin{split}
  &\matJk{\kappa}\frac{\mr{d}\uk}{\mr{d} t}+\frac{1}{2}\sum\limits_{l,m=1}^{3}
  a_{m}\left\{\Dxil{l}\matAlmk{l}{m}{\kappa}+\matAlmk{l}{m}{\kappa}\Dxil{l}\right\}\uk
  \\&-\frac{1}{2}\sum\limits_{l,m=1}^{3}\left\{
    a_{m}\diag\left(\uk\right)\Dxil{l}\matAlmk{l}{m}{\kappa}\ones{\kappa}\right\}=\\
    &\sum\limits_{l,m,a=1}^{3}\Bm{m}\Dxil{l}\matJk{\kappa}^{-1}\matAlmk{l}{m}{\kappa}\matAlmk{a}{m}{\kappa}\Dxil{a}\uk
    +\bm{SAT}_{\kappa},  
  \end{split}
\end{equation}
where $\ones{\kappa}$ is a vector of ones of size $\Nl{\kappa}$.

As in the continuous case, the semi-discrete form has a set of discrete GCL conditions
\begin{equation}\label{eq:discGCLconvection}
\sum\limits_{l=1}^{3}
    \Dxil{l}\matAlmk{l}{m}{\kappa}\ones{\kappa}=\bm{0}, \quad m = 1,2,3,
\end{equation}
that, when satisfied, lead to the following telescoping, provably stable, semi-discrete form
\begin{equation}\label{eq:convectionsplitdisctele}
  \begin{split}
  &\matJk{\kappa}\frac{\mr{d}\uk}{\mr{d} t}+\frac{1}{2}\sum\limits_{l,m=1}^{3}
  a_{m}\left\{\Dxil{l}\matAlmk{l}{m}{\kappa}+\matAlmk{l}{m}{\kappa}\Dxil{l}\right\}\uk=\\
    &\sum\limits_{l,m,a=1}^{3}\Bm{m}\Dxil{l}{\matJk{\kappa}}^{-1}\matAlmk{l}{m}{\kappa}\matAlmk{a}{m}{\kappa}\Dxil{a}\uk
    +\bm{SAT}_{\kappa}. 
  \end{split}
\end{equation}
\begin{remark}
  The linear stability of semi-discrete operators for constant coefficient hyperbolic systems, 
  is not preserved by arbitrary design order approximations of the metric terms.  Only approximations to the metric terms that satisfy
   the discrete GCL conditions~\eqref{eq:discGCLconvection} lead to stable semi-discrete forms. 
\end{remark}

\begin{remark}
The discrete metrics constructed using 
the analytic formalism of Vinokur and Yee~\cite{Vinokur2002a} or Thomas and Lombard~\cite{thomas_gcl_1979} 
will in general satisfy the discrete GCL conditions given by \eqref{eq:discGCLconvection} for
conforming interfaces when tensor-product differentiation operators are used.
\end{remark}

Herein, we optimize the metric terms as presented in \cite{fernandez_entropy_stable_p_ref_nasa_2019,fernandez_entropy_stable_p_euler_2019,fernandez_entropy_stable_p_ns_2019,fernandez_entropy_stable_hp_ref_snpdea_2019} using the 
algorithm of Crean et al. \cite{crean_entropy_stable_sbp_curvilinear_euler}. 
Only conforming interfaces are considered.
In what follows, we show which metric terms are 
optimized and how this optimization process is performed.

\subsection{Review of the inviscid coupling procedure}
The key element of the discretization that clearly show the opportunity of optimizing the 
volume metric terms is the inviscid SAT term. Thus, we consider the 
discretization of the pure convection equation, i.e., discretization
\eqref{eq:convectionsplitdisctele} with only the convective contributions:
\begin{equation}\label{eq:convectionsplitdisctele_2}
  \begin{split}
  &\matJk{\kappa}\frac{\mr{d}\uk}{\mr{d} t}+\frac{1}{2}\sum\limits_{l,m=1}^{3}
  a_{m}\left\{\Dxil{l}\matAlmk{l}{m}{\kappa}+\matAlmk{l}{m}{\kappa}\Dxil{l}\right\}\uk=
    \bm{SAT}^{(I)}_{\kappa}.
  \end{split}
\end{equation}
 \begin{figure}[t!]
   \centering
   \includegraphics[width=0.6\textwidth]{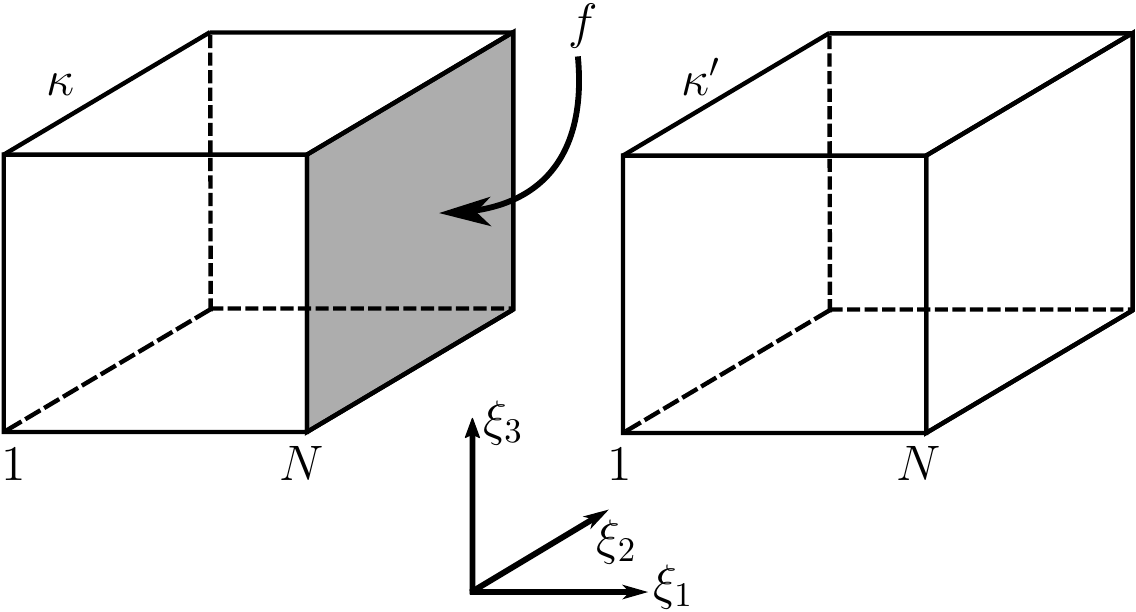}
   \caption{Generic vertical interface between two cells in computational space.}
   \label{fig:non}
 \end{figure}
Without loss of generality, and for clarity of
presentation, we consider 
only one vertical shared face $f$, as shown in Figure~\ref{fig:non}.
When considering a single set of metrics computed, for instance, with one of the approaches presented in \cite{thomas_gcl_1979} or \cite{Vinokur2002a}, 
the inviscid SAT coupling the element $\kappa$ with the neighboring element $\kappa'$ is constructed as 
\begin{equation}\label{eq:traditionalsat}
 \begin{split}
   \bm{SAT}^{(I)}_{\kappa,f} = &+ \phantom{\frac{1}{2}}{\M}^{-1}\sum\limits_{m=1}^{3}
    a_{m}\left\{\left(\eNl{1}\eNl{1}\Tr\otimes\PxiloneD{2}\otimes\PxiloneD{3}\right)\matAlmk{l}{m}{\kappa}
  \right\}\uk\\
    &-\frac{1}{2}{\M}^{-1}\sum\limits_{m=1}^{3}
    a_{m}\left\{\left(\eNl{1}\eNl{1}\Tr\otimes\PxiloneD{2}\otimes\PxiloneD{3}\right)\matAlmk{l}{m}{\kappa}
  \right\}\uk\\
    &-\frac{1}{2}{\M}^{-1}\sum\limits_{m=1}^{3}
    a_{m}\left\{\left(\eNl{1}\eonel{1}\Tr\otimes\PxiloneD{2}\otimes\PxiloneD{3}\right)\matAlmk{l}{m}{\kappa'}
  \right\}\ukp,
   \end{split}
\end{equation}
where $\M \equiv \PxiloneD{1} \otimes \PxiloneD{2} \otimes \PxiloneD{3}$.
The action of the operator $\left(\eNl{1}\eNl{1}\Tr\otimes\PxiloneD{2}\otimes\PxiloneD{3}\right)$ consists in extracting from the volume metrics
only the ``surface metrics'' associated with the LGL points located at the interface of the two cells. Similarly, $\left(\eNl{1}\eonel{1}\Tr\otimes\PxiloneD{2}\otimes\PxiloneD{3}\right)$
extracts the shared surface metrics terms from the $\kappa'$ element.
Expression \eqref{eq:traditionalsat} simplifies to the well-known formula \cite{parsani_entropy_stable_interfaces_2015}
\begin{equation}\label{eq:traditionalsat_simplified}
 \begin{split}
   \bm{SAT}^{(I)}_{\kappa,f} = &+ \frac{1}{2}{\M}^{-1}\sum\limits_{m=1}^{3}
    a_{m}\left\{\left(\eNl{1}\eNl{1}\Tr\otimes\PxiloneD{2}\otimes\PxiloneD{3}\right)\matAlmk{l}{m}{\kappa}
  \right\}\uk\\
    &-\frac{1}{2}{\M}^{-1}\sum\limits_{m=1}^{3}
    a_{m}\left\{\left(\eNl{1}\eonel{1}\Tr\otimes\PxiloneD{2}\otimes\PxiloneD{3}\right)\matAlmk{l}{m}{\kappa'}
  \right\}\ukp.
   \end{split}
\end{equation}

In this work, the convective term $\bm{SAT}^{(I)}_{\kappa,f}$
is instead constructed by using two different set of metrics
\begin{equation}
 \begin{split}
   \bm{SAT}^{(I)}_{\kappa,f} = &+ \phantom{\frac{1}{2}}{\M}^{-1}\sum\limits_{m=1}^{3}
    a_{m}\left\{\left(\eNl{1}\eNl{1}\Tr\otimes\PxiloneD{2}\otimes\PxiloneD{3}\right)\matAlmk{l}{m}{\kappa}
  \right\}\uk\\
    &-\frac{1}{2}{\M}^{-1}\sum\limits_{m=1}^{3}
    a_{m}\left\{\left(\eNl{1}\eNl{1}\Tr\otimes\PxiloneD{2}\otimes\PxiloneD{3}\right)\matAlmkAnal{l}{m}{\kappa}
  \right\}\uk\\
    &-\frac{1}{2}{\M}^{-1}\sum\limits_{m=1}^{3}
    a_{m}\left\{\left(\eNl{1}\eonel{1}\Tr\otimes\PxiloneD{2}\otimes\PxiloneD{3}\right)\matAlmkAnal{l}{m}{\kappa'}
  \right\}\ukp,
   \end{split}
\end{equation}
where the red metric terms are the analytic metrics, i.e., the metrics
computed by differentiating the inverse of the mapping \eqref{eq:3d_mapping}.
Thus, the discretization of the $\kappa$ element with the inviscid SAT contribution 
coming only from the shared vertical interface reads
\begin{equation}\label{eq:convectionsplitdisctele_1}
  \begin{split}
    \matJk{\kappa}\frac{\mr{d}\uk}{\mr{d} t} &+ \frac{1}{2}\sum\limits_{l,m=1}^{3}
  a_{m}\left\{\Dxil{l}\matAlmk{l}{m}{\kappa}+\matAlmk{l}{m}{\kappa}\Dxil{l}\right\}\uk= \\
    &+ \phantom{\frac{1}{2}}{\M}^{-1}\sum\limits_{m=1}^{3}
    a_{m}\left\{\left(\eNl{1}\eNl{1}\Tr\otimes\PxiloneD{2}\otimes\PxiloneD{3}\right)\matAlmk{l}{m}{\kappa}
  \right\}\uk\\
    &-\frac{1}{2}{\M}^{-1}\sum\limits_{m=1}^{3}
    a_{m}\left\{\left(\eNl{1}\eNl{1}\Tr\otimes\PxiloneD{2}\otimes\PxiloneD{3}\right)\matAlmkAnal{l}{m}{\kappa}
  \right\}\uk\\
    &-\frac{1}{2}{\M}^{-1}\sum\limits_{m=1}^{3}
    a_{m}\left\{\left(\eNl{1}\eonel{1}\Tr\otimes\PxiloneD{2}\otimes\PxiloneD{3}\right)\matAlmkAnal{l}{m}{\kappa'}
  \right\}\ukp.
  \end{split}
\end{equation}

Further details on how to compute the blue metric terms will be presented in the Section \ref{sec:mortar:metrics2}.

\subsection{Review of the coupling procedure for the convection-diffusion equation}\label{sec:diffusion}
In contrast to the metrics for the inviscid terms, the metrics used for the viscous terms need only be, at worst, 
consistent and design order approximations. Herein, we use the analytic metrics for the viscous terms calculation. 
To make the presentation easier and to introduce the discretization that 
will later be used for the viscous portion of the compressible Navier--Stokes equations, the inviscid term is lumped into $\fnc{I}^{(I)}$ 
while the viscous terms are simplified. Thus, Equation \eqref{eq:convectiondiffusionsplit} reduces to 
\begin{equation}\label{eq:diffusionsplit}
  \begin{split}
	  &\Jk\frac{\partial\U}{\partial t}+\fnc{I}^{(I)}= \sum\limits_{l,a=1}^{3}
  \frac{\partial}{\partial\xil{l}} \left(\Chatla{l}{a} \Thetaa{a}\right),  \\
  &\Chatla{l}{a}\equiv\sum\limits_{m=1}^{3}\Jk\frac{\partial\xil{l}}{\partial\xm{m}}\frac{\partial\xil{a}}{\partial\xm{m}}\Bm{m},\quad 
  \Thetaa{a}\equiv\frac{\partial\U}{\partial\xil{a}}.
  \end{split}
\end{equation}
 A local discontinuous Galerkin (LDG) and interior penalty approach (IP) approach are used (see references~\cite{parsani_entropy_stable_interfaces_2015,carpenter_entropy_stable_staggered_2015,parsani_ssdc_staggered_2016}). 
 In the LDG approach, the discretization of the viscous terms in Equation \eqref{eq:diffusionsplit} proceeds in two steps. First, the gradients $\Thetaa{a}$ are discretized, 
 then the derivatives of the viscous fluxes are discretized. Notice that all the metric terms are contained in $\Chatla{l}{a}$, and therefore the critical 
 requirement for stability is to use an SBP operator \cite{parsani_entropy_stable_interfaces_2015,carpenter_entropy_stable_staggered_2015,fernandez_entropy_stable_p_ref_nasa_2019}. Plugging everything together, the final discretization reads
 \begin{equation}\label{eq:diffusionsplitdisc}
	 \matJk{\kappa}\frac{\mr{d}\uk}{\mr{d}t}+\bm{I}^{(I)}_{\kappa}=\sum\limits_{l,a=1}^{3}\Dxil{l}\matChatlaAnal{l}{a}_{\kappa}\thetaa{a}^{\kappa}\:+\: \bm{SAT}^{(I)}_{\kappa} + \bm{SAT}^{(V)}_{\kappa} ,\quad
	 \thetaa{a}^{\kappa}=\Dxil{a}\uk + \bm{SAT}^{\theta}_{\kappa},
 \end{equation}
 where the inviscid contributions are contained in $\bm{I}^{(I)}_{\kappa}$, while $\bm{SAT}^{\theta}_{\kappa}$ contains the LDG penalty on the gradient of the entropy variables \cite{parsani_entropy_stable_interfaces_2015}.
The viscous coefficient matrices $\matChatlaAnal{l}{a}_{\kappa}$ have been highlighted in red to emphasize that they are computed using the analytic metrics. However, any other design order choice would suffice.
The proposed discretization of the viscous terms telescopes the 
 viscous fluxes to the boundary and adds a dissipative term \cite{parsani_entropy_stable_interfaces_2015}. Thus, it mimics the continuous energy analysis, and leads to
 a provably energy stable discretization, provided appropriate boundary SATs are used.



\subsection{Metric solution mechanics}\label{sec:mortar:metrics2}
Here we demonstrate the proposed approach for the approximation of the metric 
terms, and we consider the discrete GCL conditions \eqref{eq:discGCLconvection} associated with Equation
\eqref{eq:convectionsplitdisctele_1}.
By construction, the analytic metric terms are continuous at the interface of the two elements
because the curvilinear coordinate transformation \eqref{eq:3d_mapping} is water-tight (see Assumption \ref{assume:curv}).
Thus the GCL constraints read 
\begin{equation}\label{eq:discGCLconvection_2}
  \begin{split}
    \M\sum\limits_{l=1}^{3}\Dxil{l}\matAlmk{l}{m}{\kappa}\ones{\kappa}=&
   + 
    \left\{\left(\eNl{1}\eNl{1}\Tr\otimes\PxiloneD{2}\otimes\PxiloneD{3}\right)\matAlmk{l}{m}{\kappa}
  \right\}\ones{\kappa}\\
    &-
    \left\{\left(\eNl{1}\eonel{1}\Tr\otimes\PxiloneD{2}\otimes\PxiloneD{3}\right)\matAlmkAnal{l}{m}{\kappa'}
  \right\}\ones{\kappa}, \\
   &m = 1,2,3.
  \end{split}
\end{equation}

Equation \eqref{eq:discGCLconvection_2} can be algebraically manipulated into a 
form that is more convenient for constructing a solution procedure for the metric terms. 
Multiplying Equation \eqref{eq:discGCLconvection_2} by $-1$, using the SBP property $\QxiloneD{l} = -{\QxiloneD{l}}\Tr+\ExiloneD{l}$,
and by canceling common terms, we arrive at
\begin{equation}\label{eq:discGCLconvection_3}
  \begin{split}
    \sum\limits_{l=1}^{3}\Qxil{l}\Tr\matAlmk{l}{m}{\kappa}\ones{\kappa}=
    \left\{\left(\eNl{1}\eonel{1}\Tr\otimes\PxiloneD{2}\otimes\PxiloneD{3}\right)\matAlmkAnal{l}{m}{\kappa'}
  \right\}\ones{\kappa},\quad 
   m = 1,2,3,
  \end{split}
\end{equation}
where 
$\Qxil{1}\equiv\QxiloneD{1}\otimes\PxiloneD{2}\otimes\PxiloneD{3}$, 
$\Qxil{2}\equiv\PxiloneD{1}\otimes\QxiloneD{2}\otimes\PxiloneD{3}$,
$\Qxil{3}\equiv\PxiloneD{1}\otimes\PxiloneD{2}\otimes\QxiloneD{3}$.

A close examination of \eqref{eq:discGCLconvection_3}
shows that the GCL constraints form a highly under-determined system for the unknown metric terms on the left-hand side. 
Our strategy consists of solving
a strictly convex quadratic optimization problem that minimizes the difference between the numerical and analytic volume metrics \cite{crean_entropy_stable_sbp_curvilinear_euler}:

\begin{equation}\label{eq:opt}
  \begin{split}
&\min_{\bm{a}_{m}}\frac{1}{2}\left(\bm{a}_{m}
-\bm{a}_{m,\text{target}}\right)\Tr\left(\bm{a}_{m}-\bm{a}_{m,\text{target}}\right),\quad
\text{subject to}\quad\mat{M}\bm{a}_{m}=\bm{c}_{m},\\
 &m = 1,2,3,
  \end{split}
\end{equation}
where $\bm{a}_{m}$ is a vector of size $3\Nl{\kappa}$ containing the optimized metrics, \ie,
\begin{equation*}
\left(\bm{a}_{m}\right)\Tr\equiv\ones{}\Tr
\left[
\matAlmk{1}{m}{\kappa},
\matAlmk{2}{m}{\kappa},
\matAlmk{3}{m}{\kappa}
\right]
\end{equation*}
and $\bm{a}_{m,\text{target}}$ are the targeted, analytical metrics, \ie
\begin{equation*}
\left(\bm{a}_{m,\text{target}}\right)\Tr\equiv\ones{}\Tr
\left[
\matAlmkAnal{1}{m}{\kappa},
\matAlmkAnal{2}{m}{\kappa},
\matAlmkAnal{3}{m}{\kappa}
\right].
\end{equation*}

The constraints $\mat{M}\bm{a}_{m}=\bm{c}_{m}$ are simply 
the discrete GCL conditions \eqref{eq:discGCLconvection_3}. Specifically, the matrix $\mat{M}$ is of size $\Nl{\kappa}\times 3 \Nl{\kappa}$,
and it is defined as
\begin{equation}
\mat{M} \equiv\left[\Qxil{1}\Tr,\Qxil{2}\Tr,\Qxil{3}\Tr\right],
\end{equation}
while the right-hand side data for the constrained equations, $\bm{c}_{m}$, is a vector of size 
 $\Nl{\kappa}$ defined as
\begin{equation}\label{eq:cm}
  \bm{c}_{m} \equiv \frac{1}{2}\left\{\left(\eNl{1}\eonel{1}\Tr\otimes\PxiloneD{2}\otimes\PxiloneD{3}\right)\matAlmkAnal{l}{m}{\kappa'}
  \right\}\ones{\kappa}.
\end{equation}
The optimal solution of the constrained minimization problem 
is given by (see Proposition $1$ in \cite{crean_entropy_stable_sbp_curvilinear_euler})
\begin{equation}\label{eq:opt1}
\bm{a}_{m}^{}=\bm{a}_{m,\text{target}}^{}-{\MM^{}}^{\dagger}\left(\MM^{} 
\bm{a}_{m,\text{target}}^{}-\bm{c}_{m}^{}\right),
\end{equation}
with ${\MM^{}}^{\dagger}$ the Moore-Penrose pseudo-inverse of $\MM^{}$. 
\section{Discretization of the compressible Navier--Stokes equations}\label{sec:compressible_nse}
In this section, the algorithm 
for the convection-diffusion equation presented in the previous section is applied 
to the compressible Navier--Stokes equations with conforming interfaces. 
These equations in Cartesian coordinates read 
\begin{equation}\label{eq:compressible_ns}
\begin{aligned}
  &\frac{\partial\Q}{\partial t}+\sum\limits_{m=1}^{3}\frac{\partial \FxmI{m}}{\partial \xm{m}} = \sum\limits_{m=1}^{3}\frac{\partial \FxmV{m}}{\partial\xm{m}}, &&
\forall \left(\xm{1},\xm{2},\xm{3}\right)\in\Omega,\quad t\ge 0,\\
  &\Q\left(\xm{1},\xm{2},\xm{3},t\right)=\GB\left(\xm{1},\xm{2},\xm{3},t\right), && \forall \left(\xm{1},\xm{2},\xm{3}\right)\in\Gamma,\quad t\ge 0,\\
&\Q\left(\xm{1},\xm{2},\xm{3},0\right)=\Gzero\left(\xm{1},\xm{2},\xm{3}\right), &&
  \forall \left(\xm{1},\xm{2},\xm{3}\right)\in\Omega,
\end{aligned}
\end{equation}
where the vectors $\Q$, $\FxmI{m}$ and $\FxmV{m}$ denote the conserved variables, the inviscid fluxes, and the
viscous fluxes, respectively. The boundary data,
$\GB$, and the initial condition, $\Gzero$, are assumed to be in $L^{2}(\Omega)$,
with the further assumption that $\GB$ will be set to coincide with linear,
well-posed boundary conditions, prescribed in such a way that either entropy conservation or entropy stability is achieved.

The vector of conserved variables is given by 
\begin{equation*}
\Q = \left[\rho,\rho\Um{1},\rho\Um{2},\rho\Um{3},\rho\E\right]\Tr,
\end{equation*}
where $\rho$ denotes the density, $\bm{\fnc{U}} = \left[\Um{1},\Um{2},\Um{3}\right]\Tr$ is the velocity 
vector, and $\E$ is the specific total energy. The inviscid fluxes are given as
\begin{equation*}
  \begin{split}
\FxmI{m} = &\left[\rho\Um{m},\rho\Um{m}\Um{1}+\delta_{m,1}\fnc{P},\rho\Um{m}\Um{2}+\delta_{m,2}\fnc{P},\right.
\left.\rho\Um{m}\Um{3}+\delta_{m,3}\fnc{P},\rho\Um{m}\fnc{H}\right]\Tr,
  \end{split}
\end{equation*}
where $\fnc{P}$ is the pressure, $\fnc{H}$ is the specific total enthalpy and $\delta_{i,j}$ is the 
Kronecker delta.

The required constituent relations are
\begin{equation*}
\fnc{H} = c_{\fnc{P}}\fnc{T}+\frac{1}{2}\bm{\fnc{U}}\Tr\bm{\fnc{U}},\quad \fnc{P} = \rho R \fnc{T},\quad R = \frac{R_{u}}{M_{w}},
\end{equation*}
where $\fnc{T}$ is the temperature, $R_{u}$ is the universal gas constant, $M_{w}$ is the molecular weight of the gas, 
and $c_{\fnc{P}}$ is the specific heat capacity at constant pressure. Finally, the specific thermodynamic entropy is given as 
\begin{equation*}
s=\frac{R}{\gamma-1}\log\left(\frac{\fnc{T}}{\fnc{T}_{\infty}}\right)-R\log\left(\frac{\rho}{\rho_{\infty}}\right),\quad \gamma=\frac{c_{p}}{c_{p}-R},
\end{equation*}
where $\fnc{T}_{\infty}$ and $\rho_{\infty}$ are the reference temperature and density, respectively
(the stipulated convention has been broken here and $s$ has been used rather than $\fnc{S}$ for reasons that will be clear next). 

The viscous fluxes $\FxmV{m}$ are given by
\begin{equation}\label{eq:Fv}
\FxmV{m}=\left[0,\tau_{1,m},\tau_{2,m},\tau_{3,m},
\sum\limits_{i=1}^{3}\tau_{i,m}\fnc{U}_{i}-\kappa\frac{\partial \fnc{T}}{\partial\xm{m}}\right]\Tr,
\end{equation} 
while the viscous stresses are defined as
\begin{equation}\label{eq:tau}
\tau_{i,j} = \mu\left(\frac{\partial\fnc{U}_{i}}{\partial x_{j}}+\frac{\partial\fnc{U}_{j}}{\partial x_{i}}
-\delta_{i,j}\frac{2}{3}\sum\limits_{n=1}^{3}\frac{\partial\fnc{U}_{n}}{\partial x_{n}}\right),
\end{equation}
where $\mu(\fnc{T})$ is the dynamic viscosity and $\kappa(\fnc{T})$ is the thermal conductivity.

The compressible Navier--Stokes equations given in \eqref{eq:NSCCS1} have
a convex extension, that when integrated over the physical domain, $\Omega$, 
depends only on the boundary data and negative semi-definite dissipation terms.
This convex extension depends on an entropy function, $\fnc{S}$, 
that is constructed from the thermodynamic entropy as
\[
\fnc{S}=-\rho s,
\]
and provides a mechanism for proving stability in the $L^{2}$ norm.  
The entropy variables $\bfnc{W}$ are an alternative variable set related to the conservative variables via a one-to-one mapping.  They
are defined in terms of the entropy function $\fnc{S}$ by the relation
$\bfnc{W}\Tr=\partial\fnc{S}/\partial\bfnc{Q}$ and they are extensively used in the
entropy stability proofs of the algorithms used herein; see for instance
\cite{carpenter_ssdc_2014,parsani_ssdc_staggered_2016,friedrich_hp_entropy_stability_2018,fernandez_entropy_stable_hp_ref_snpdea_2019}.
In addition, they simultaneously symmetrize
the inviscid and the viscous flux Jacobians in all three spatial directions.
Further details on continuous entropy analysis are available 
elsewhere \cite{dafermos_book_2010,parsani_entropy_stability_solid_wall_2015,carpenter_entropy_stable_staggered_2015}.

The entropy stability for the viscous terms in the compressible Navier--Stokes 
equations \eqref{eq:NSCCS1} is readily demonstrated by 
exploiting the symmetrizing properties of the 
entropy variables.
Thus, we recast the 
viscous fluxes in terms of the entropy variables
\begin{equation}\label{eq:Fxment}
\FxmV=\sum\limits_{j=1}^{3}\Cij{m}{j}\frac{\partial\bfnc{W}}{\partial x_{j}},
\end{equation}
with the flux Jacobian matrices satisfying $\Cij{m}{j} \:=\: {(\Cij{j}{m})}\Tr$.

Furthermore, in order to apply the algorithm 
outlined for the convection-diffusion case \eqref{eq:diffusionsplitdisc} to the compressible Navier--Stokes equations, we have
to recast system \eqref{eq:compressible_ns} in
a skew-symmetric form with respect to the metric terms. This procedure results in
\begin{equation}\label{eq:NSCCS1}
\begin{split}
&\Jk\frac{\partial\bfnc{Q}}{\partial t}+\sum\limits_{l,m=1}^{3}\frac{1}{2}
\frac{\partial }{\partial \xil{l}}\left(\Jdxildxm{l}{m}\FxmI{m}\right)
+\frac{1}{2}\Jdxildxm{l}{m}\frac{\partial \FxmI{m}}{\partial \xil{l}}
=\sum\limits_{l,m=1}^{3}\frac{\partial}{\partial\xil{l}}\left(\Jdxildxm{l}{m}\FxmV{m}\right)
\end{split},
\end{equation}
where the GCL relations given in \eqref{eq:GCL} are used to obtain \eqref{eq:NSCCS} from the divergence form \eqref{eq:compressible_ns}.
Substituting \eqref{eq:Fxment} into~\eqref{eq:NSCCS1},  we arrive at the system of equations
\begin{equation}\label{eq:NSCCS}
\begin{split}
&\Jk\frac{\partial\bfnc{Q}}{\partial t}+\sum\limits_{l,m=1}^{3}
\frac{1}{2}\frac{\partial }{\partial \xil{l}}\left(\Jdxildxm{l}{m}\Fxm{l}\right)
+\frac{1}{2}\Jdxildxm{l}{m}\frac{\partial \Fxm{m}}{\partial \xil{l}}
=\sum\limits_{l,a=1}^{3}\frac{\partial}{\partial\xil{l}}\left(\Chatij{l}{a}\frac{\partial\bfnc{W}}{\partial \xil{a}}\right),
\end{split}
\end{equation}
where
\begin{equation}\label{eq:Chatij}
\Chatij{l}{a}=\Jdxildxm{l}{m}\sum\limits_{m,j=1}^{3}\Cij{m}{j}\frac{\partial\xil{a}}{\partial x_{j}}.
\end{equation}
The symmetric properties of the viscous flux Jacobians are preserved by the rotation into curvilinear
coordinates, \ie $\Chatij{l}{a} \:=\: {(\Chatij{a}{l})}\Tr$.
We remark that this form of the equations, \ie skew-symmetric form plus the quadratic form of the viscous terms, 
is necessary for the construction of the entropy stable schemes used in this work.
For further details on the derivation of these viscous coefficient matrices see \cite{fisher_phd_2012,parsani_entropy_stability_solid_wall_2015}.

Without loss of generality, as was done for the linear convection-diffusion equation, 
we consider only the coupling SAT terms for one shared interface. Thus, the 
discretization of the compressible Euler equations, i.e., the inviscid part of \eqref{eq:NSCCS}, is given by 
\begin{equation}\label{eq:discEL}
  \begin{split}
    &\matJk{\kappa}\frac{\mr{d}\qk{\kappa}}{\mr{d}t}+\frac{1}{2}\sum\limits_{l,m=1}^{3}\left(\Dxil{l}\matAlmk{l}{m}{\kappa}+\matAlmk{m}{l}{\kappa}\Dxil{l}\right)
  \circ\matFxm{m}{\qk{\kappa}}{\qk{\kappa}}\ones{\kappa}=\\
    &+\phantom{\frac{1}{2}}\M^{-1}\sum\limits_{m=1}^{3}
    \left\{\left(\eNl{1}\eNl{1}\Tr\otimes\PxiloneD{2}\otimes\PxiloneD{3}\otimes\Imat{5}\right)\matAlmk{l}{m}{\kappa}
  \right\} \circ\matFxm{m}{\qk{\kappa}}{\qk{\kappa}}\ones{\kappa}\\
    &- \frac{1}{2}\M^{-1}\sum\limits_{m=1}^{3}
    \left\{\left(\eNl{1}\eNl{1}\Tr\otimes\PxiloneD{2}\otimes\PxiloneD{3}\otimes\Imat{5}\right)\matAlmkAnal{l}{m}{\kappa}
  \right\} \circ\matFxm{m}{\qk{\kappa}}{\qk{\kappa'}}\ones{\kappa}\\
    &-\frac{1}{2}\M^{-1}\sum\limits_{m=1}^{3}
    \left\{\left(\eNl{1}\eonel{1}\Tr\otimes\PxiloneD{2}\otimes\PxiloneD{3}\otimes\Imat{5}\right)\matAlmkAnal{l}{m}{\kappa'}
    \right\} \circ\matFxm{m}{\qk{\kappa}}{\qk{\kappa'}}\ones{\kappa},
  \end{split}
\end{equation}
where the symbol $\circ$ indicates the Hadamard product, and $\matFxm{m}{\cdot}{\cdot}$ 
is a two argument matrix flux function which is constructed from a two point entropy
conservative flux function (see, for instance, \cite{fernandez_entropy_stable_p_ref_nasa_2019}).
The Hadamard formalism is capable of compactly representing various split forms, 
and more importantly, extends to nonlinear equations for which a canonical split form is inappropriate.
It is used in the construction of entropy conservative/stable discretizations which 
are used herein.

Next, recasting the viscous fluxes in terms of entropy variables as shown in \eqref{eq:Fxment} yields the following form for the discretization of the divergence of the viscous fluxes
\begin{equation}\label{eq:MacroForm}
\sum\limits_{l,a=1}^{3}\frac{\partial}{\partial\xil{l}}\left(\Chatla{l}{a}\frac{\partial\bfnc{W}}{\partial \xil{a}}\right)
	\approx  \sum\limits_{l,a=1}^{3} \Dxil{l}\matChatlaAnal{l}{a}\thetaa{a}^{\kappa}, \qquad \thetaa{a}^{\kappa}=\Dxil{a}\wk{\kappa}. 
\end{equation}
Note that Equation \eqref{eq:MacroForm} is precisely the symmetric generalization of the convection-diffusion operator 
to a viscous system.  

The discretization on the $\kappa\Th$ element reads
\begin{equation}\label{eq:NSL}
  \begin{split}
	  &\matJk{\kappa}\frac{\mr{d}\qk{\kappa}}{\mr{d}t}+\bm{I}^{(E)}_{\kappa}=\sum\limits_{l,a=1}^{3}
	  \Dxil{a}\matChatlaAnal{l}{a}\thetaa{a}^{\kappa} + \bm{SAT}^{(I)}_{\kappa} + \bm{SAT}^{(V)}_{\kappa} + \IP^{\kappa},\quad
	  \thetaa{a}^{\kappa}=\Dxil{a}\wk{\kappa} + \bm{SAT}^{\theta}_{\kappa},
  \end{split}
\end{equation}
where $\bm{I}^{(E)}_{\kappa}$ represents the discretization of the divergence of the inviscid fluxes and the interior penalty term, $\IP^{\kappa}$, adds interface dissipation \cite{parsani_entropy_stable_interfaces_2015}.
This term is a design-order zero interface dissipation term that is constructed to damp neutrally stable ``odd-even'' eigenmodes that arise from the LDG viscous operator. Scheme \eqref{eq:NSL} telescopes to the boundaries where appropriate SATs need to be imposed to 
obtain a stability statement \cite{parsani_entropy_stable_interfaces_2015,carpenter_entropy_stable_staggered_2015,parsani_ssdc_staggered_2016}.

\section{Numerical results}\label{sec:numerical_results} 
Herein, the conforming \cite{carpenter_ssdc_2014,parsani_wall_bc_entropy_2015,carpenter_entropy_stability_ssdc_2016,parsani_ssdc_staggered_2016} and $p$-adaptive solver
\cite{fernandez_entropy_stable_p_ref_nasa_2019,fernandez_entropy_stable_p_ns_2019,fernandez_entropy_stable_p_euler_2019} for unstructured grids 
developed at the Extreme Computing Research Center (ECRC) at KAUST is used to perform numerical experiments. This parallel solver is built on 
top of the Portable and Extensible Toolkit for Scientific computing (PETSc)~\cite{petsc-user-ref}, its mesh topology 
abstraction (DMPLEX)~\cite{KnepleyKarpeev09} and scalable ordinary differential equation (ODE)/differential algebraic equations (DAE) solver library~\cite{abhyankar2018petsc}. 
The systems of ordinary differential equations arising from the spatial
discretizations are integrated using the fourth-order
accurate Bogacki--Shampine method \cite{BOGACKI1989321} endowed with an adaptive time stepping technique based on digital signal processing \cite{Soderlind2003,Soderlind2006}. To make the temporal error negligible, a tolerance of $10^{-8}$ is always used for the time-step adaptivity. The two-point entropy consistent flux
of Chandrashekar~\cite{Chandrashekar2013} is used for all the test cases.

The errors are computed using a volume scaled $L^{2}$ discrete norm as follows:
\begin{equation*}
\begin{split}
  &\|\bm{u}\|_{L^{2}}^{2}=\Omega_{c}^{-1}\sum\limits_{\kappa=1}^{K}\bm{u}_{\kappa}\M\matJk{\kappa}\bm{u}_{\kappa},
\end{split}
\end{equation*}
where $\Omega_{c}$ indicates the volume of $\Omega$ computed as 
$\Omega_{c}\equiv\sum\limits_{\kappa=1}^{K}\ones{\kappa}\Tr\M\matJk{\kappa}\ones{\kappa}$.

We study the $L^2$ norm of the error in the primitive variables (density, 
velocity components, and temperature) considering two test cases having \textit{analytical} 
solution, and using two different computational domains
discretized with geometrically high order grids.

\subsection{Isentropic vortex}\label{subsec:isentripic_vortex}

In this section, we report on the numerical results for the
propagation of an isentropic vortex by solving the three-dimensional compressible
Euler equations. The analytical solution of this problem is
\begin{equation} \label{vortex-solution}
\begin{split}
& \fnc{G} = 1
-\left\{
\left[
  \left(\xm{1}-x_{1,0}\right)
-U_{\infty}\cos\left(\alpha\right)t
\right]^{2}
+
\left[
  \left(\xm{2}-x_{2,0}\right)
-U_{\infty}\sin\left(\alpha\right)t
\right]^{2}
\right\},\\
&\rho = \fnc{T}^{\frac{1}{\gamma-1}},
 \quad
 \fnc{T} = \left[1-\epsilon_{\nu}^{2}M_{\infty}^{2}\frac{\gamma-1}{8\pi^{2}}\exp\left(\fnc{G}\right)\right],\\
  &\Um{1} = U_{\infty}\cos(\alpha)-\epsilon_{\nu}
  \frac{\left(\xm{2}-x_{2,0}\right)-U_{\infty}\sin\left(\alpha\right)t}{2\pi}
\exp\left(\frac{\fnc{G}}{2}\right),\\
  &\U{2} = U_{\infty}\sin(\alpha)-\epsilon_{\nu}
  \frac{\left(\xm{1}-x_{1,0}\right)-U_{\infty}\cos\left(\alpha\right)t}{2\pi}
\exp\left(\frac{\fnc{G}}{2}\right), \\
	&\Um{3} = 1,
\end{split}
\end{equation}
where $U_{\infty}$ is the modulus of the free-stream velocity, $M_{\infty}$
is the free-stream Mach number, $c_{\infty}$ is the free-stream speed of sound,
and $\left(x_{1,0},x_{2,0},x_{3,0}\right)$ is the vortex center.
The following values are used: $U_{\infty}=M_{\infty} c_{\infty}$, $\epsilon_{\nu}=5$, $M_{\infty}=0.5$,
$\gamma=1.4$, $\alpha=45^{\degree}$, and $\left(x_{1,0},x_{2,0},x_{3,0}\right)=\left(0,0,0\right)$.
The initial condition is given by \eqref{vortex-solution} with $t=0$.

\subsubsection{Cubic domain with geometrically high order perturbed cells}\label{subsec:iv_cube}

\begin{figure}
     \centering
     \subfloat[][Sample of the grid.]{\includegraphics[width=.39\textwidth]{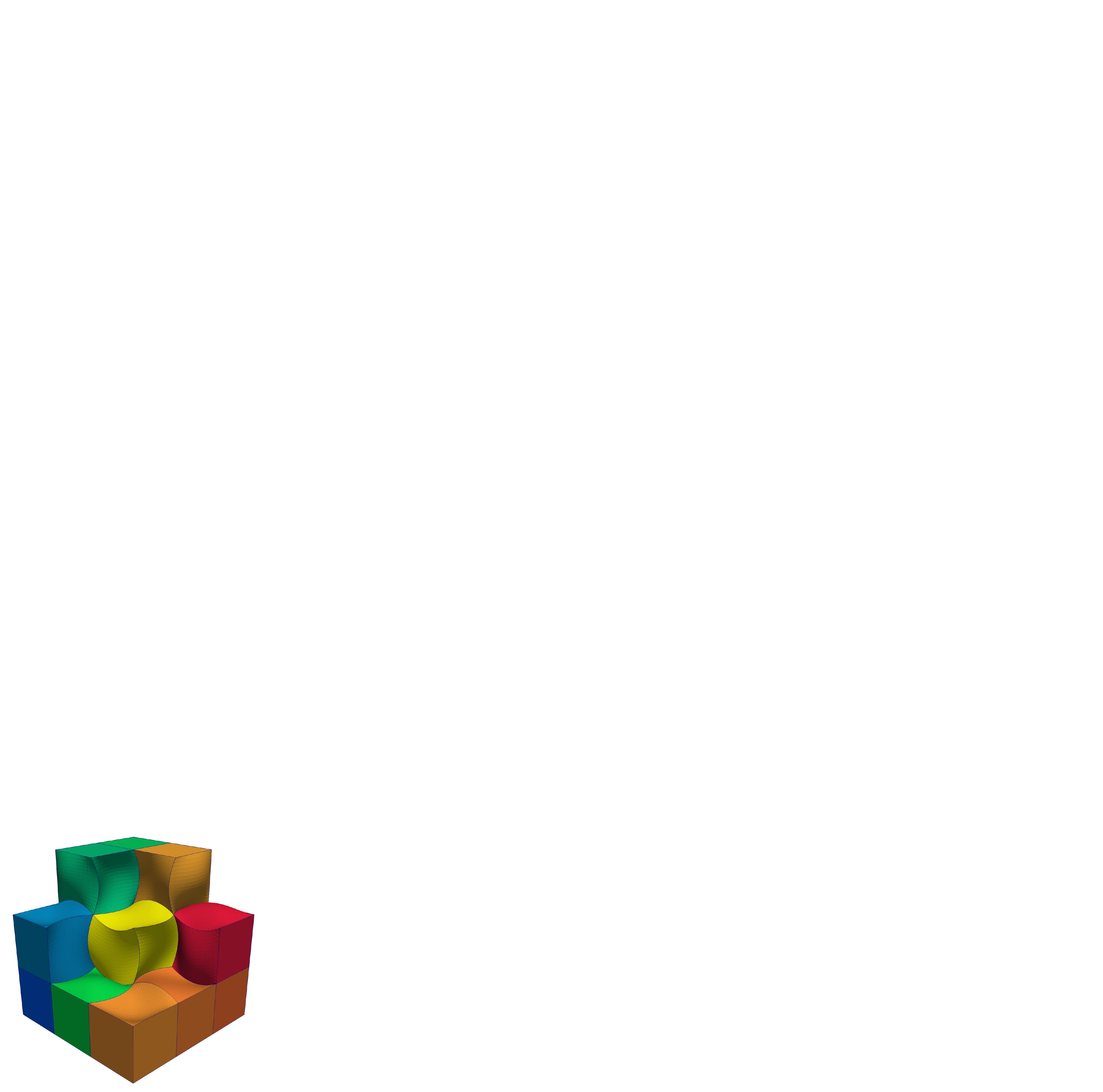}\label{fig:sample_perturb}}
     \subfloat[][Density of the isentropic vortex at $t=0$.]{\includegraphics[width=.55\textwidth]{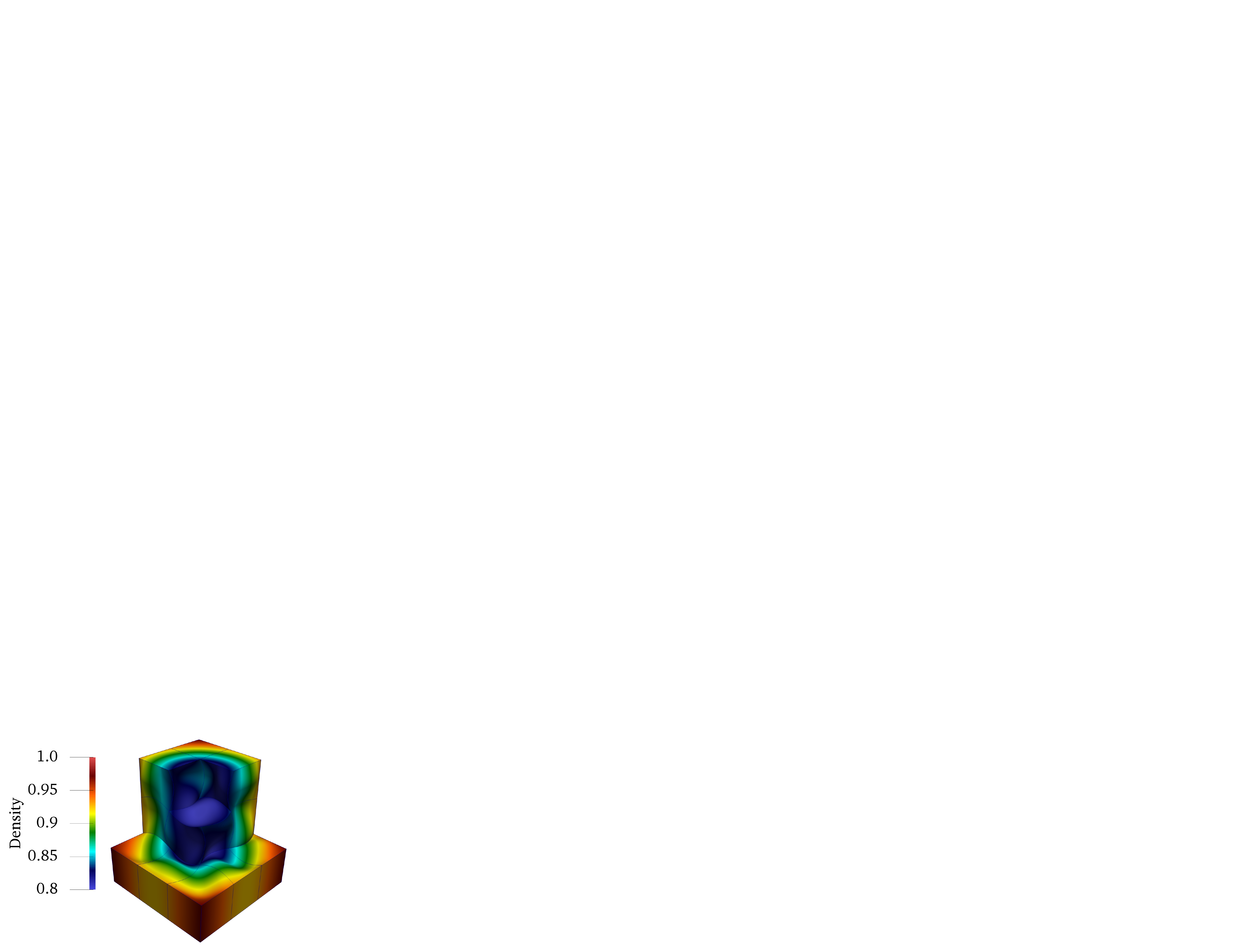}\label{fig:v_mag_iv}}
     \caption{Cubic domain with geometrically high order perturbed cells.}
     \label{fig:cube_perturb}
\end{figure}

The physical domain is $[-1,1]^3$, uniformly subdivided into 27 hexahedrons.
We first collocate the LGL points in physical coordinates by an affine mapping of the reference element.
These physical nodal coordinates, denoted by $x_{i,*}$, are then perturbed as
\begin{equation*}
\begin{split}
&x_1 = x_{1,*} + 2 \frac{\eta}{15} \cos \left(   a \right) \cos \left( 3 b \right) \sin \left( 4 c \right),\\
&x_2 = x_{2,*} + 2 \frac{\eta}{15} \sin \left( 4 a \right) \cos \left(   b \right) \cos \left( 3 c \right),\\
&x_3 = x_{3,*} + 2 \frac{\eta}{15} \cos \left( 3 a \right) \sin \left( 4 b \right) \cos \left(   c \right),
\end{split}
\end{equation*}
where $a = \frac{\pi}{2} x_{1,*}$, $b = \frac{\pi}{2} x_{2,*}$, and $c = \frac{\pi}{2} x_{3,*}$,
while $\eta$ is a positive perturbation parameter $0\leq\eta\leq 1$. Figure \ref{fig:sample_perturb} shows a cut of the mesh for $\eta=1$.

In Table \ref{tbl:ivortex_cube}, we report the ratio of the $L^2$ norm of the 
errors of the primitive variables computed using the metrics of Thomas and Lombard  \cite{thomas_gcl_1979} and the 
optimization algorithm presented herein. Solution polynomials of degree $p=1,2,3,4,7,15$ 
and four values of the perturbation parameter, i.e., $\eta = 0.25, 0.5, 0.75, 1$, are considered. 
The numbers highlighted in green indicate that the 
ratio of the error norm is equal to or greater than one, i.e., for that specific 
primitive variable the solution computed with the optimized metric terms is more accurate than that
obtained by using the metrics of Thomas and Lombard \cite{thomas_gcl_1979}. 

We
observe that $95\%$ of the ratios are equal to or greater than one with some cases where 
the error is reduced by a factor greater than $2.5$. For the cases 
where the ratio is highlighted in red, the maximum deviation from one is approximately $0.05\%$.
Furthermore, we note that
the ratio of the two error norms converges to one when increasing the degree of the solution polynomial.

\begin{table}
\centering
  \subfloat[Density.]{
\begin{tabular}{l|c|c|c|c}
 & \multicolumn{4}{c}{$\rho$}  \\ 
\cline{1-5}
$\eta$          & $0.25$ & $0.5$ & $0.75$  & $1.0$  \\ 
\hline
$p=1 $ & \cellcolor{green!25} 1.000 &\cellcolor{green!25} 1.000 & \cellcolor{green!25} 1.000 & \cellcolor{green!25} 1.000  \\
$p=2 $ & \cellcolor{green!25} 1.291 &\cellcolor{green!25} 1.557 & \cellcolor{green!25} 1.606 & \cellcolor{green!25} 1.534  \\
$p=3 $ & \cellcolor{green!25} 2.821 &\cellcolor{green!25} 2.570 & \cellcolor{green!25} 2.137 & \cellcolor{green!25} 1.796  \\
$p=4 $ & \cellcolor{green!25} 2.916 &\cellcolor{green!25} 2.343 & \cellcolor{green!25} 1.882 & \cellcolor{green!25} 1.582  \\
$p=7 $ & \cellcolor{green!25} 1.346 &\cellcolor{green!25} 1.111 & \cellcolor{green!25} 1.051 & \cellcolor{green!25} 1.026  \\ 
$p=9 $ & \cellcolor{green!25} 1.047 &\cellcolor{green!25} 1.013 & \cellcolor{green!25} 1.003 & \cellcolor{red!25}   0.998  \\
$p=15$ & \cellcolor{green!25} 1.015 &\cellcolor{green!25} 1.000 & \cellcolor{green!25} 1.000 & \cellcolor{green!25} 1.000           
\end{tabular}
}
\qquad
  \subfloat[Velocity component in $\xm{1}$.]{
\begin{tabular}{l|c|c|c|c}
 & \multicolumn{4}{c}{$\mathcal{U}_{1}$}  \\ 
\cline{1-5}
$\eta$          & $0.25$ & $0.5$ & $0.75$  & $1.0$  \\ 
\hline
$p=1 $ &\cellcolor{green!25} 1.000 & \cellcolor{green!25} 1.000 & \cellcolor{green!25} 1.000 & \cellcolor{green!25} 1.000 \\
$p=2 $ &\cellcolor{green!25} 1.485 & \cellcolor{green!25} 1.849 & \cellcolor{green!25} 1.835 & \cellcolor{green!25} 1.721 \\
$p=3 $ &\cellcolor{green!25} 1.829 & \cellcolor{green!25} 1.775 & \cellcolor{green!25} 1.578 & \cellcolor{green!25} 1.444 \\
$p=4 $ &\cellcolor{green!25} 1.583 & \cellcolor{green!25} 1.319 & \cellcolor{green!25} 1.207 & \cellcolor{green!25} 1.162 \\
$p=7 $ &\cellcolor{green!25} 1.044 & \cellcolor{green!25} 1.034 & \cellcolor{green!25} 1.031 & \cellcolor{green!25} 1.019 \\ 
$p=9 $ &\cellcolor{green!25} 1.009 & \cellcolor{green!25} 1.006 & \cellcolor{green!25} 1.001 & \cellcolor{red!25}   0.996 \\ 
$p=15$ &\cellcolor{green!25} 1.001 & \cellcolor{green!25} 1.000 & \cellcolor{green!25} 1.000 & \cellcolor{green!25} 1.000        
\end{tabular}
}
\qquad
\subfloat[Velocity component in $\xm{2}$.]{
\begin{tabular}{l|c|c|c|c}
 & \multicolumn{4}{c}{$\mathcal{U}_{2}$}  \\ 
\cline{1-5}
$\eta$          & $0.25$ & $0.5$ & $0.75$  & $1.0$  \\
\hline
$p=1 $ &\cellcolor{green!25} 1.000 &\cellcolor{green!25} 1.000 &\cellcolor{green!25} 1.000 &\cellcolor{green!25} 1.000 \\
$p=2 $ &\cellcolor{green!25} 1.543 &\cellcolor{green!25} 2.027 &\cellcolor{green!25} 2.078 &\cellcolor{green!25} 1.966 \\
$p=3 $ &\cellcolor{green!25} 1.822 &\cellcolor{green!25} 1.807 &\cellcolor{green!25} 1.602 &\cellcolor{green!25} 1.436 \\
$p=4 $ &\cellcolor{green!25} 1.974 &\cellcolor{green!25} 1.556 &\cellcolor{green!25} 1.333 &\cellcolor{green!25} 1.232 \\
$p=7 $ &\cellcolor{green!25} 1.056 &\cellcolor{green!25} 1.038 &\cellcolor{green!25} 1.025 &\cellcolor{green!25} 1.017 \\ 
$p=9 $ &\cellcolor{green!25} 1.017 &\cellcolor{green!25} 1.005 &\cellcolor{red!25}   0.998 &\cellcolor{red!25}   0.995 \\
$p=15$ &\cellcolor{red!25}   0.999 &\cellcolor{red!25}   0.999 &\cellcolor{green!25} 1.000 &\cellcolor{green!25} 1.000            
\end{tabular}
}
\qquad
\subfloat[Velocity component in $\xm{3}$.]{
\begin{tabular}{l|c|c|c|c}
 & \multicolumn{4}{c}{$\mathcal{U}_{3}$}  \\ 
\cline{1-5}
$\eta$          & $0.25$ & $0.5$ & $0.75$  & $1.0$  \\ 
\hline
$p=1 $ &\cellcolor{green!25} 1.000 &\cellcolor{green!25} 1.000 &\cellcolor{green!25} 1.000 &\cellcolor{green!25} 1.000  \\
$p=2 $ &\cellcolor{green!25} 1.365 &\cellcolor{green!25} 1.420 &\cellcolor{green!25} 1.271 &\cellcolor{green!25} 1.130  \\
$p=3 $ &\cellcolor{green!25} 2.228 &\cellcolor{green!25} 1.850 &\cellcolor{green!25} 1.560 &\cellcolor{green!25} 1.391  \\
$p=4 $ &\cellcolor{green!25} 2.212 &\cellcolor{green!25} 1.682 &\cellcolor{green!25} 1.380 &\cellcolor{green!25} 1.217  \\
$p=7 $ &\cellcolor{green!25} 1.194 &\cellcolor{green!25} 1.074 &\cellcolor{green!25} 1.043 &\cellcolor{green!25} 1.026  \\ 
$p=9 $ &\cellcolor{green!25} 1.032 &\cellcolor{green!25} 1.010 &\cellcolor{green!25} 1.001 &\cellcolor{red!25}   0.997  \\
$p=15$ &\cellcolor{green!25} 1.000 &\cellcolor{green!25} 1.000 &\cellcolor{green!25} 1.000 &\cellcolor{green!25} 1.000         
\end{tabular}
}
\qquad
\subfloat[Temperature.]{
\begin{tabular}{l|c|c|c|c}
 & \multicolumn{4}{c}{$\mathcal{T}$}  \\ 
\cline{1-5}
$\eta$          & $0.25$ & $0.5$ & $0.75$  & $1.0$  \\ 
\hline
$p=1 $ & \cellcolor{green!25} 1.000 & \cellcolor{green!25} 1.000 & \cellcolor{green!25} 1.000 & \cellcolor{green!25} 1.000  \\
$p=2 $ & \cellcolor{green!25} 1.237 & \cellcolor{green!25} 1.473 & \cellcolor{green!25} 1.547 & \cellcolor{green!25} 1.528  \\
$p=3 $ & \cellcolor{green!25} 2.364 & \cellcolor{green!25} 2.434 & \cellcolor{green!25} 2.150 & \cellcolor{green!25} 1.893  \\
$p=4 $ & \cellcolor{green!25} 2.281 & \cellcolor{green!25} 1.989 & \cellcolor{green!25} 1.691 & \cellcolor{green!25} 1.478  \\
$p=7 $ & \cellcolor{green!25} 1.299 & \cellcolor{green!25} 1.098 & \cellcolor{green!25} 1.047 & \cellcolor{green!25} 1.026  \\ 
$p=9 $ & \cellcolor{green!25} 1.033 & \cellcolor{green!25} 1.010 & \cellcolor{green!25} 1.003 & \cellcolor{green!25} 1.000  \\
$p=15$ & \cellcolor{green!25} 1.005 & \cellcolor{green!25} 1.000 & \cellcolor{green!25} 1.000 & \cellcolor{green!25} 1.000         
\end{tabular}
}
\caption{Ratio of the $L^2$-norm of the errors (Thomas and Lombard \cite{thomas_gcl_1979} vs. optimized): 
  propagation of the isentropic vortex in a cubic domain with geometrically high order perturbed cells.}
\label{tbl:ivortex_cube}
\end{table}

\subsubsection{Cubic domain with a spherical hole and geometrically high-order cells}\label{subsec:iv_sphere}

\begin{figure}
 \begin{center}
   \includegraphics[width=0.5\textwidth]{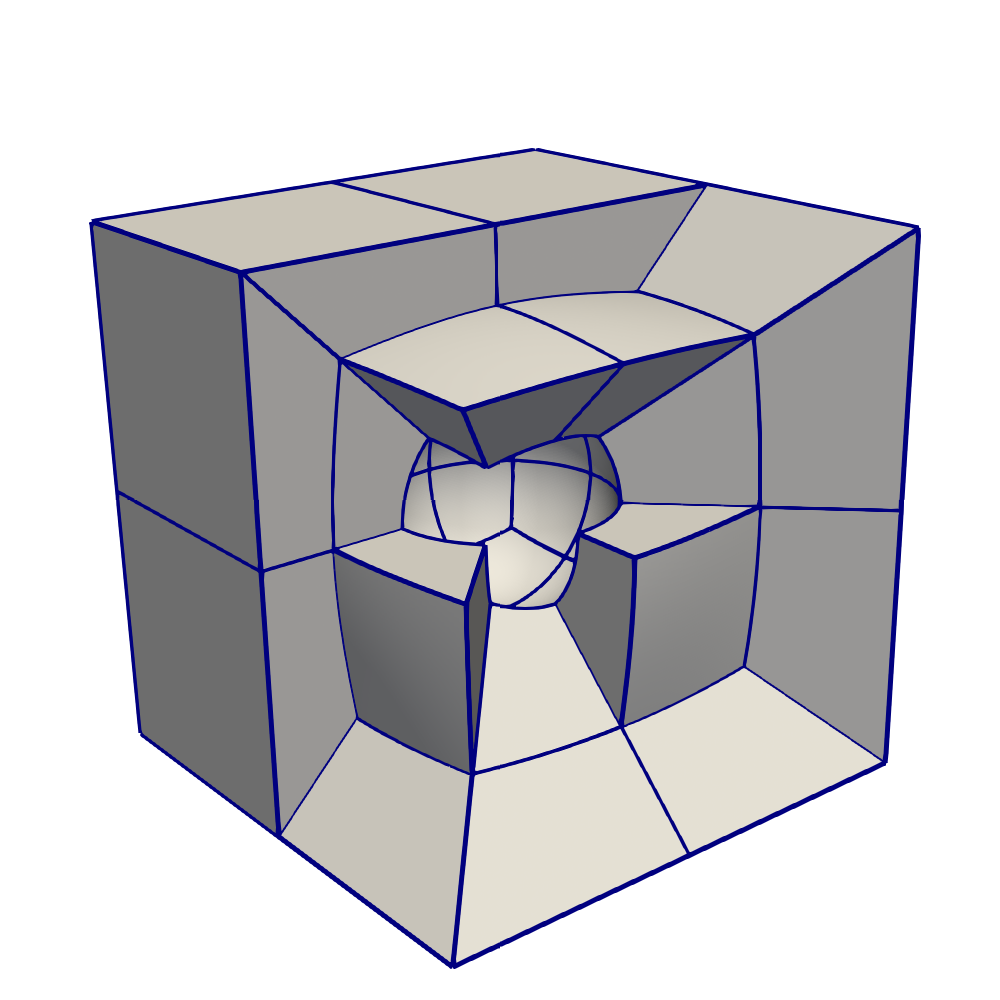}
   \caption{Cubic domain with a spherical hole and curved cells}
   \label{fig:box_hole}
  \end{center}
\end{figure}

In this section, we present the results for the propagation of the
isentropic vortex in a cubic domain $[-1,1]^3$ with a spherical hole
of radius $0.5$ placed at the center. The domain is approximated with
a mesh consisting of $48$ hexahedrons; a cut of the mesh is detailed in Figure \ref{fig:box_hole}.

In Table \ref{tbl:ivortex_cube_sphere}, we report the ratio of the $L^2$ norm of the 
errors at a final time $t_{f}=1$ using solution polynomials of degree $p=2,3,4,7,9$.
Both interior and boundary faces describing the spherical geometry are approximated
with polynomials of the same order.
We observe that $92\%$ of the ratios are equal to or greater than one. In this case,
where the mesh is less distorted, the maximum improvement is approximately $10\%$.
The ratios highlighted in red deviate from one by at most $0.01\%$.

\begin{table}[htbp!]
\begin{center}
\begin{tabular}{l|c|c|c|c|c}
  &$\rho$          & $\mathcal{U}_{1}$ & $\mathcal{U}_{2}$ & $\mathcal{U}_{3}$  & $\mathcal{T}$  \\ 
\hline
 $ p=2  $ & \cellcolor{green!25} 1.094  & \cellcolor{red!25}    0.999   & \cellcolor{red!25}   0.999  & \cellcolor{green!25} 1.139  & \cellcolor{green!25} 1.093\\
 $ p=3  $ & \cellcolor{green!25} 1.018  & \cellcolor{green!25}  1.019   & \cellcolor{green!25} 1.019  & \cellcolor{green!25} 1.007  & \cellcolor{green!25} 1.024\\
 $ p=4  $ & \cellcolor{green!25} 1.057  & \cellcolor{green!25}  1.002   & \cellcolor{green!25} 1.002  & \cellcolor{green!25} 1.028  & \cellcolor{green!25} 1.016\\
 $ p=7  $ & \cellcolor{green!25} 1.041  & \cellcolor{green!25}  1.014   & \cellcolor{green!25} 1.014  & \cellcolor{green!25} 1.039  & \cellcolor{green!25} 1.023\\
 $ p=9  $ & \cellcolor{green!25} 1.029  & \cellcolor{green!25}  1.011   & \cellcolor{green!25} 1.011  & \cellcolor{green!25} 1.017  & \cellcolor{green!25} 1.022          
\end{tabular}
\caption{Ratio of the $L^2$-norm of the errors (Thomas and Lombard \cite{thomas_gcl_1979} vs. optimized): 
  propagation of the isentropic vortex in a cubic domain with a spherical hole and geometrically high order cells.}        
\label{tbl:ivortex_cube_sphere}
\end{center}
\end{table}

\subsection{Viscous shock}\label{subsec:viscous_shock}
Next, we study the propagation of a viscous shock using the
compressible Navier--Stokes equations.  We assume a planar shock propagating along the $\xm{1}$
coordinate direction with a Prandtl number of $Pr=3/4$.
The exact solution of this problem is known;
the momentum $\fnc{V}(x_1)$ satisfies the ordinary differential equation
\begin{equation*}
\begin{split}
  &\alpha\fnc{V}\frac{\partial\fnc{V}}{\partial\xm{1}}-(\fnc{V}-1)(\fnc{V}-\fnc{V}_{f})=0, \qquad -\infty\leq\xm{1}\leq+\infty.
\end{split}
\end{equation*}
Assuming that the center of the viscous shock is located at $\xm{1} = 0$, the implicit solution of the former equation is
\begin{equation}\label{eq:implicit_sol_vs}
  \xm{1}-\frac{1}{2}\alpha\left(\log\left|(\fnc{V}(x_1)-1)(\fnc{V}(x_1)-\fnc{V}_{f})\right|+\frac{1+\fnc{V}_{f}}{1-\fnc{V}_{f}}\log\left|\frac{\fnc{V}(x_1)-1}{\fnc{V}(x_1)-\fnc{V}_{f}}\right|\right) = 0,
\end{equation}
where
\begin{equation}
\fnc{V}_{f}\equiv\frac{\fnc{U}_{L}}{\fnc{U}_{R}},\qquad
\alpha\equiv\frac{2\gamma}{\gamma + 1}\frac{\,\mu}{Pr\dot{\fnc{M}}}.
\end{equation}
Here $\fnc{U}_{L/R}$ are known velocities to the left and right of the shock at
$-\infty$ and $+\infty$, respectively, $\dot{\fnc{M}}$ is the constant mass
flow across the shock, $Pr$ is the Prandtl number, and $\mu$ is the dynamic
viscosity.
The mass and total enthalpy are constant across the shock, and
the momentum and energy equations become redundant.

For our tests, $\fnc{V}$ is computed from Equation \eqref{eq:implicit_sol_vs}
to machine precision using bisection.
The moving shock solution is obtained by applying a uniform translation to the above solution.
The shock is located at the center of the domain at $t=0$, and the following
values are used: $M_{\infty}=2.5$, $Re_{\infty}=10$, and $\gamma=1.4$. The boundary conditions are prescribed by penalizing the numerical solution against the exact solution. The analytical solution is also used to furnish data for the initial conditions.

We investigate the $L^2$ norm of the error for the same meshes discussed for the isentropic vortex configurations; the simulations are stopped at a final time $t_{f} = 5$.

\subsubsection{Cubic domain with geometrically high order perturbed cells}\label{subsec:vs_cube}
We consider solution polynomials of degree $p=1,2,3,4,7,15$ 
and four values of the perturbation parameter, i.e., $\eta = 0.25, 0.5, 0.75, 1$. In Table \ref{tbl:vs_cube}, we report the ratio of the $L^2$ norm of the 
errors of the primitive variables computed using the metrics of Thomas and Lombard  \cite{thomas_gcl_1979} and the 
optimization algorithm presented herein. We
observe that $80\%$ of the ratios are equal to or greater than one, where in some cases the improvement is more than $40\%$. For the cases 
where the ratio is highlighted in red, the maximum deviation from one is approximately $1.6\%$.
Furthermore, we note that by increasing the degree of the solution polynomial
the ratio of the two error norms converges to one.

\begin{table}[htbp!]
\centering
  \subfloat[]{
\begin{tabular}{l|c|c|c|c}
 & \multicolumn{4}{c}{$\rho$}  \\ 
\cline{1-5}
$\eta$          & $0.25$ & $0.5$ & $0.75$  & $1.0$  \\ 
\hline
$p=1 $ & \cellcolor{green!25} 1.000 & \cellcolor{green!25} 1.000 & \cellcolor{green!25} 1.000 & \cellcolor{green!25} 1.000 \\
$p=2 $ & \cellcolor{green!25} 1.213 & \cellcolor{green!25} 1.536 & \cellcolor{green!25} 1.776 & \cellcolor{green!25} 1.951 \\
$p=3 $ & \cellcolor{green!25} 1.027 & \cellcolor{green!25} 1.055 & \cellcolor{green!25} 1.078 & \cellcolor{green!25} 1.098 \\
$p=4 $ & \cellcolor{green!25} 1.050 & \cellcolor{green!25} 1.072 & \cellcolor{green!25} 1.058 & \cellcolor{green!25} 1.032 \\
$p=7 $ & \cellcolor{red!25}   0.996 & \cellcolor{red!25}   0.995 & \cellcolor{red!25}   0.996 & \cellcolor{red!25}   0.999 \\ 
$p=9 $ & \cellcolor{green!25} 1.000 & \cellcolor{green!25} 1.000 & \cellcolor{green!25} 1.000 & \cellcolor{green!25} 1.000 \\
$p=15$ &\cellcolor{green!25} 1.000 & \cellcolor{green!25} 1.000 & \cellcolor{green!25} 1.000 & \cellcolor{green!25} 1.000  
\end{tabular}
}
\qquad
\subfloat[]{
\begin{tabular}{l|c|c|c|c}
 & \multicolumn{4}{c}{$\mathcal{U}_{1}$}  \\ 
\cline{1-5}
$\eta$          & $0.25$ & $0.5$ & $0.75$  & $1.0$  \\ 
\hline
$p=1 $ & \cellcolor{green!25}  1.000 & \cellcolor{green!25} 1.000 & \cellcolor{green!25} 1.000 & \cellcolor{green!25} 1.000 \\
$p=2 $ & \cellcolor{green!25}  1.187 & \cellcolor{green!25} 1.418 & \cellcolor{green!25} 1.546 & \cellcolor{green!25} 1.623 \\
$p=3 $ & \cellcolor{green!25}  1.034 & \cellcolor{green!25} 1.054 & \cellcolor{green!25} 1.059 & \cellcolor{green!25} 1.067 \\
$p=4 $ & \cellcolor{green!25}  1.098 & \cellcolor{green!25} 1.121 & \cellcolor{green!25} 1.114 & \cellcolor{green!25} 1.098 \\
$p=7 $ & \cellcolor{red!25}    0.992 & \cellcolor{red!25}   0.986 & \cellcolor{red!25}   0.984 & \cellcolor{red!25}   0.985 \\
$p=9 $ & \cellcolor{red!25}    0.998 & \cellcolor{red!25}   0.998 & \cellcolor{red!25}   0.998 & \cellcolor{red!25}   0.999 \\
$p=15$ &\cellcolor{green!25}  1.000 & \cellcolor{green!25} 1.000 & \cellcolor{green!25} 1.000 & \cellcolor{green!25} 1.000
\end{tabular}
}
\qquad
\subfloat[]{
\begin{tabular}{l|c|c|c|c}
 & \multicolumn{4}{c}{$\mathcal{U}_{2}$}  \\ 
\cline{1-5}
$\eta$          & $0.25$ & $0.5$ & $0.75$  & $1.0$  \\
\hline
$p=1 $ & \cellcolor{green!25} 1.000 & \cellcolor{green!25} 1.000 & \cellcolor{green!25} 1.000 & \cellcolor{green!25} 1.000 \\
$p=2 $ & \cellcolor{green!25} 1.187 & \cellcolor{green!25} 1.418 & \cellcolor{green!25} 1.546 & \cellcolor{green!25} 1.623 \\
$p=3 $ & \cellcolor{green!25} 1.034 & \cellcolor{green!25} 1.054 & \cellcolor{green!25} 1.059 & \cellcolor{green!25} 1.067 \\
$p=4 $ & \cellcolor{green!25} 1.098 & \cellcolor{green!25} 1.121 & \cellcolor{green!25} 1.114 & \cellcolor{green!25} 1.098 \\
$p=7 $ & \cellcolor{red!25}   0.992 & \cellcolor{red!25}   0.986 & \cellcolor{red!25}   0.984 & \cellcolor{red!25}   0.985 \\
$p=9 $ & \cellcolor{red!25}   0.998 & \cellcolor{red!25}   0.998 & \cellcolor{red!25}   0.998 & \cellcolor{red!25}   0.999 \\
$p=15$ &\cellcolor{green!25} 1.000 & \cellcolor{green!25} 1.000 & \cellcolor{green!25} 1.000 & \cellcolor{green!25} 1.000       
\end{tabular}
}
\qquad
\subfloat[]{
\begin{tabular}{l|c|c|c|c}
 & \multicolumn{4}{c}{$\mathcal{U}_{3}$}  \\ 
\cline{1-5}
$\eta$          & $0.25$ & $0.5$ & $0.75$  & $1.0$  \\ 
\hline
$p=1 $ & \cellcolor{green!25} 1.000 & \cellcolor{green!25} 1.000 & \cellcolor{green!25} 1.000 & \cellcolor{green!25} 1.000 \\
$p=2 $ & \cellcolor{green!25} 1.187 & \cellcolor{green!25} 1.418 & \cellcolor{green!25} 1.546 & \cellcolor{green!25} 1.623 \\
$p=3 $ & \cellcolor{green!25} 1.034 & \cellcolor{green!25} 1.054 & \cellcolor{green!25} 1.059 & \cellcolor{green!25} 1.067 \\
$p=4 $ & \cellcolor{green!25} 1.098 & \cellcolor{green!25} 1.121 & \cellcolor{green!25} 1.114 & \cellcolor{green!25} 1.098 \\
$p=7 $ & \cellcolor{red!25}   0.992 & \cellcolor{red!25}   0.986 & \cellcolor{red!25}   0.984 & \cellcolor{red!25}   0.985 \\ 
$p=9 $ & \cellcolor{red!25}   0.998 & \cellcolor{red!25}   0.998 & \cellcolor{red!25}   0.998 & \cellcolor{red!25}   0.999 \\
$p=15$ &\cellcolor{green!25} 1.000 & \cellcolor{green!25} 1.000 & \cellcolor{green!25} 1.000 & \cellcolor{green!25} 1.000    
\end{tabular}
}
\qquad
\subfloat[]{
\begin{tabular}{l|c|c|c|c}
 & \multicolumn{4}{c}{$\mathcal{T}$}  \\ 
\cline{1-5}
$\eta$          & $0.25$ & $0.5$ & $0.75$  & $1.0$  \\ 
\hline
$p=1 $ & \cellcolor{green!25} 1.000 &  \cellcolor{green!25} 1.000 & \cellcolor{green!25} 1.000 & \cellcolor{green!25} 1.000 \\
$p=2 $ & \cellcolor{green!25} 1.227 &  \cellcolor{green!25} 1.467 & \cellcolor{green!25} 1.568 & \cellcolor{green!25} 1.633 \\
$p=3 $ & \cellcolor{green!25} 1.105 &  \cellcolor{green!25} 1.182 & \cellcolor{green!25} 1.226 & \cellcolor{green!25} 1.262 \\
$p=4 $ & \cellcolor{green!25} 1.203 &  \cellcolor{green!25} 1.208 & \cellcolor{green!25} 1.184 & \cellcolor{green!25} 1.161 \\
$p=7 $ & \cellcolor{red!25}   0.996 &  \cellcolor{green!25} 1.011 & \cellcolor{green!25} 1.020 & \cellcolor{green!25} 1.023 \\ 
$p=9 $ & \cellcolor{green!25} 1.002 &  \cellcolor{green!25} 1.006 & \cellcolor{green!25} 1.006 & \cellcolor{green!25} 1.005 \\
$p=15$ &\cellcolor{green!25} 1.000 &  \cellcolor{green!25} 1.000 & \cellcolor{green!25} 1.000 & \cellcolor{green!25} 1.000         
\end{tabular}
}
\caption{Ratio of the $L^2$-norm of the errors (Thomas and Lombard \cite{thomas_gcl_1979} vs. optimized): 
  propagation of the viscous shock in a cubic domain with geometrically high order perturbed cells.}
\label{tbl:vs_cube}
\end{table}

\subsubsection{Cubic domain with a spherical hole and geometrically high order cells}\label{subsec:vs_sphere}
In Table \ref{tbl:vs_cube_sphere}, we report the ratio of the $L^2$ norm of the 
errors of the primitive variables. Solution polynomials of degree $p=2,3,4,7,9$ 
are used. Boundary faces describing the spherical geometry as well as internal curved 
interfaces are approximated
with polynomials of degree $p$.
We
observe that $70\%$ of the ratios are equal to or greater than one with a maximum gain
of $5\%$ for the solution computed with the optimized metrics.
The ratio highlighted in red has a maximum deviation from one of just $1.9\%$.

\begin{table}[htbp!]
\begin{center}
\begin{tabular}{l|c|c|c|c|c}
&$\rho$          & $\mathcal{U}_{1}$ & $\mathcal{U}_{2}$ & $\mathcal{U}_{3}$  & $\mathcal{T}$  \\ 
\hline
$p=2 $ & \cellcolor{green!25} 1.013    & \cellcolor{red!25}   0.981   &  \cellcolor{red!25}   0.981  & \cellcolor{red!25}   0.981     &\cellcolor{green!25} 1.047\\
$p=3 $ & \cellcolor{green!25} 1.000    & \cellcolor{red!25}   0.997   &  \cellcolor{red!25}   0.997  & \cellcolor{red!25}   0.997     &\cellcolor{red!25}   0.996\\
$p=4 $ & \cellcolor{green!25} 1.001    & \cellcolor{green!25} 1.000   &  \cellcolor{green!25} 1.000  & \cellcolor{green!25} 1.000     &\cellcolor{green!25} 1.005\\
$p=7 $ & \cellcolor{green!25} 1.000    & \cellcolor{green!25} 1.000   &  \cellcolor{green!25} 1.000  & \cellcolor{green!25} 1.000     &\cellcolor{green!25} 1.000\\
$p=9 $ & \cellcolor{green!25} 1.000    & \cellcolor{green!25} 1.000   &  \cellcolor{green!25} 1.000  & \cellcolor{green!25} 1.000     &\cellcolor{green!25} 1.000           
\end{tabular}
\caption{Ratio of the $L^2$ norm of the errors (Thomas and Lombard \cite{thomas_gcl_1979} vs. optimized): 
  propagation of the viscous shock in a cubic domain with a spherical hole and geometrically high order cells.}        
\label{tbl:vs_cube_sphere}
\end{center}
\end{table}

\section{Conclusion}\label{sec:conclusions}
We numerically show that optimizing the metric terms as proposed in 
\cite{fernandez_entropy_stable_p_ref_nasa_2019,fernandez_entropy_stable_p_euler_2019,fernandez_entropy_stable_hp_ref_snpdea_2019}
lead, even for conforming interfaces, to a solution whose accuracy is 
practically never worse and often noticeably better than the one obtained using 
the widely adopted Thomas and Lombard metric terms computation \cite{thomas_gcl_1979}.
We also observed that by increasing the degree of the solution polynomial,
the ratio of the two error norms converges to one. This indicates that the benefits
of optimizing the metric terms decrease by increasing the solution polynomial degree of
the spatial approximation.
We conclude that the pre-processing step of optimizing the metric terms 
can then be used in a computational framework as a unique and viable approach for conforming and $h/p$ non-conforming
interfaces. In addition, this choice greatly simplifies the solver and allows important code re-utilization.
Investigating the effect of the optimized metrics on other important systems of PDEs is a future research direction.
\section*{Acknowledgments}
The research reported in this paper was funded by King Abdullah University of
Science and Technology. We are thankful for the computing resources of the
Supercomputing Laboratory and the Extreme Computing Research Center at
King Abdullah University of Science and Technology.

\bibliographystyle{siam}
\bibliography{metrics}

\pagebreak

\appendix

\section{$L^2$ norm of the error: Isentropic vortex}
\subsection{Cubic domain with geometrically high order perturbed cells}

\begin{table}
\centering
  \subfloat[Thomas and Lombard metrics \cite{thomas_gcl_1979}.]{
\begin{tabular}{l|c|c|c|c}
 & \multicolumn{4}{c}{$\rho$}  \\ 
\cline{1-5}
$\eta$          & $0.25$ & $0.5$ & $0.75$  & $1.0$  \\ 
\hline
$p=1 $ & 6.062e-03 & 6.062e-03 & 6.062e-03 & 6.062e-03 \\
$p=2 $ & 1.210e-03 & 1.888e-03 & 2.610e-03 & 3.313e-03 \\
$p=3 $ & 4.024e-04 & 7.994e-04 & 1.210e-03 & 1.633e-03 \\
$p=4 $ & 1.456e-04 & 2.982e-04 & 4.683e-04 & 6.644e-04 \\
$p=7 $ & 1.381e-06 & 4.875e-06 & 1.193e-05 & 2.398e-05 \\
$p=9 $ & 6.974e-08 & 3.724e-07 & 1.112e-06 & 2.506e-06 \\
$p=15$ & 1.127e-11 & 1.186e-10 & 5.584e-10 & 1.847e-09       
\end{tabular}
}
\qquad
  \subfloat[Optimized metrics.]{
\begin{tabular}{l|c|c|c|c}
 & \multicolumn{4}{c}{$\rho$}  \\ 
\cline{1-5}
$\eta$          & $0.25$ & $0.5$ & $0.75$  & $1.0$  \\ 
\hline
$p=1 $ & 6.062e-03 & 6.062e-03 & 6.062e-03 & 6.062e-03 \\
$p=2 $ & 9.371e-04 & 1.212e-03 & 1.625e-03 & 2.161e-03 \\
$p=3 $ & 1.427e-04 & 3.110e-04 & 5.664e-04 & 9.097e-04 \\
$p=4 $ & 4.993e-05 & 1.273e-04 & 2.488e-04 & 4.199e-04 \\
$p=7 $ & 1.026e-06 & 4.388e-06 & 1.135e-05 & 2.338e-05 \\
$p=9 $ & 6.660e-08 & 3.677e-07 & 1.109e-06 & 2.511e-06 \\
$p=15$ & 1.110e-11 & 1.186e-10 & 5.585e-10 & 1.847e-09            
\end{tabular}
}
\caption{$L^2$ error norm of the density: 
  propagation of the isentropic vortex in a cubic domain with geometrically high order perturbed cells.}
\label{tbl:ivortex_cube_rho}
\end{table}

\begin{table}
\centering
  \subfloat[Thomas and Lombard metrics \cite{thomas_gcl_1979}.]{
\begin{tabular}{l|c|c|c|c}
  & \multicolumn{4}{c}{$\Um{1}$}  \\ 
\cline{1-5}
$\eta$          & $0.25$ & $0.5$ & $0.75$  & $1.0$  \\ 
\hline
$p=1 $ & 1.694e-02 & 1.694e-02 & 1.694e-02 & 1.694e-02 \\
$p=2 $ & 4.635e-03 & 7.990e-03 & 1.170e-02 & 1.562e-02 \\
$p=3 $ & 1.146e-03 & 2.393e-03 & 4.003e-03 & 6.003e-03 \\
$p=4 $ & 3.589e-04 & 8.494e-04 & 1.563e-03 & 2.534e-03 \\
$p=7 $ & 6.144e-06 & 2.445e-05 & 5.768e-05 & 1.117e-04 \\
$p=9 $ & 3.691e-07 & 1.673e-06 & 4.549e-06 & 1.022e-05 \\
$p=15$ &  2.608e-11 & 3.024e-10 & 1.685e-09 & 6.562e-09      
\end{tabular}
}
\qquad
  \subfloat[Optimized metrics.]{
\begin{tabular}{l|c|c|c|c}
  & \multicolumn{4}{c}{$\Um{1}$}  \\ 
\cline{1-5}
$\eta$          & $0.25$ & $0.5$ & $0.75$  & $1.0$  \\ 
\hline
$p=1 $ & 1.694e-02 & 1.694e-02 & 1.694e-02 & 1.694e-02  \\
$p=2 $ & 3.121e-03 & 4.322e-03 & 6.376e-03 & 9.074e-03  \\
$p=3 $ & 6.269e-04 & 1.348e-03 & 2.537e-03 & 4.156e-03  \\
$p=4 $ & 2.267e-04 & 6.441e-04 & 1.295e-03 & 2.180e-03  \\
$p=7 $ & 5.886e-06 & 2.364e-05 & 5.597e-05 & 1.096e-04  \\
$p=9 $ & 3.659e-07 & 1.663e-06 & 4.543e-06 & 1.027e-05  \\
$p=15$ &  2.606e-11 & 3.024e-10 & 1.685e-09 & 6.562e-09            
\end{tabular}
}
  \caption{$L^2$ error norm of the velocity component in $\xm{1}$: 
  propagation of the isentropic vortex in a cubic domain with geometrically high order perturbed cells.}
\label{tbl:ivortex_cube_u1}
\end{table}

\begin{table}
\centering
  \subfloat[Thomas and Lombard metrics \cite{thomas_gcl_1979}.]{
\begin{tabular}{l|c|c|c|c}
  & \multicolumn{4}{c}{$\Um{2}$}  \\ 
\cline{1-5}
$\eta$          & $0.25$ & $0.5$ & $0.75$  & $1.0$  \\ 
\hline
$p=1 $ & 1.694e-02 & 1.694e-02 & 1.694e-02 & 1.694e-02 \\
$p=2 $ & 4.731e-03 & 8.215e-03 & 1.207e-02 & 1.614e-02 \\
$p=3 $ & 1.279e-03 & 2.590e-03 & 4.183e-03 & 6.158e-03 \\
$p=4 $ & 4.159e-04 & 9.325e-04 & 1.638e-03 & 2.595e-03 \\
$p=7 $ & 6.125e-06 & 2.372e-05 & 5.962e-05 & 1.215e-04 \\
$p=9 $ & 3.406e-07 & 1.738e-06 & 5.138e-06 & 1.162e-05 \\
$p=15$ &  4.037e-11 & 4.686e-10 & 2.518e-09 & 9.149e-09      
\end{tabular}
}
\qquad
  \subfloat[Optimized metrics.]{
\begin{tabular}{l|c|c|c|c}
  & \multicolumn{4}{c}{$\Um{2}$}  \\ 
\cline{1-5}
$\eta$          & $0.25$ & $0.5$ & $0.75$  & $1.0$  \\ 
\hline
$p=1 $ & 1.694e-02 & 1.694e-02 & 1.694e-02 & 1.694e-02 \\
$p=2 $ & 3.065e-03 & 4.053e-03 & 5.807e-03 & 8.212e-03 \\
$p=3 $ & 7.021e-04 & 1.433e-03 & 2.612e-03 & 4.287e-03 \\
$p=4 $ & 2.107e-04 & 5.995e-04 & 1.229e-03 & 2.107e-03 \\
$p=7 $ & 5.799e-06 & 2.285e-05 & 5.818e-05 & 1.195e-04 \\
$p=9 $ & 3.349e-07 & 1.730e-06 & 5.148e-06 & 1.168e-05 \\
$p=15$ &  4.039e-11 & 4.689e-10 & 2.519e-09 & 9.151e-09           
\end{tabular}
}
  \caption{$L^2$ error norm of the velocity component in $\xm{2}$: 
  propagation of the isentropic vortex in a cubic domain with geometrically high order perturbed cells.}
\label{tbl:ivortex_cube_u2}
\end{table}

\begin{table}
\centering
  \subfloat[Thomas and Lombard metrics \cite{thomas_gcl_1979}.]{
\begin{tabular}{l|c|c|c|c}
  & \multicolumn{4}{c}{$\Um{3}$}  \\ 
\cline{1-5}
$\eta$          & $0.25$ & $0.5$ & $0.75$  & $1.0$  \\ 
\hline
$p=1 $ & 7.813e-03 & 7.813e-03 & 7.813e-03 & 7.813e-03 \\
$p=2 $ & 2.439e-03 & 4.468e-03 & 6.427e-03 & 8.357e-03 \\
$p=3 $ & 1.272e-03 & 2.577e-03 & 4.051e-03 & 5.800e-03 \\
$p=4 $ & 4.767e-04 & 1.009e-03 & 1.648e-03 & 2.448e-03 \\
$p=7 $ & 5.851e-06 & 2.218e-05 & 5.467e-05 & 1.102e-04 \\
$p=9 $ & 3.187e-07 & 1.746e-06 & 5.195e-06 & 1.180e-05 \\
$p=15$ &  4.049e-11 & 4.332e-10 & 2.107e-09 & 7.318e-09      
\end{tabular}
}
\qquad
  \subfloat[Optimized metrics.]{
\begin{tabular}{l|c|c|c|c}
  & \multicolumn{4}{c}{$\Um{3}$}  \\ 
\cline{1-5}
$\eta$          & $0.25$ & $0.5$ & $0.75$  & $1.0$  \\ 
\hline
$p=1 $ & 7.813e-03 & 7.813e-03 & 7.813e-03 & 7.813e-03 \\
$p=2 $ & 1.786e-03 & 3.145e-03 & 5.057e-03 & 7.393e-03 \\
$p=3 $ & 5.711e-04 & 1.393e-03 & 2.597e-03 & 4.170e-03 \\
$p=4 $ & 2.155e-04 & 6.000e-04 & 1.195e-03 & 2.012e-03 \\
$p=7 $ & 4.900e-06 & 2.065e-05 & 5.241e-05 & 1.073e-04 \\
$p=9 $ & 3.089e-07 & 1.729e-06 & 5.191e-06 & 1.183e-05 \\
$p=15$ &  4.049e-11 & 4.333e-10 & 2.107e-09 & 7.319e-09           
\end{tabular}
}
  \caption{$L^2$ error norm of the velocity component in $\xm{3}$: 
  propagation of the isentropic vortex in a cubic domain with geometrically high order perturbed cells.}
\label{tbl:ivortex_cube_u3}
\end{table}

\begin{table}
\centering
  \subfloat[Thomas and Lombard metrics \cite{thomas_gcl_1979}.]{
\begin{tabular}{l|c|c|c|c}
  & \multicolumn{4}{c}{$\mathcal{T}$}  \\ 
\cline{1-5}
$\eta$          & $0.25$ & $0.5$ & $0.75$  & $1.0$  \\ 
\hline
$p=1 $ & 3.038e-03 & 3.038e-03 & 3.038e-03 & 3.038e-03 \\
$p=2 $ & 5.380e-04 & 8.319e-04 & 1.181e-03 & 1.562e-03 \\
$p=3 $ & 1.838e-04 & 3.624e-04 & 5.555e-04 & 7.683e-04 \\
$p=4 $ & 6.521e-05 & 1.342e-04 & 2.139e-04 & 3.095e-04 \\
$p=7 $ & 6.432e-07 & 2.329e-06 & 5.680e-06 & 1.130e-05 \\
$p=9 $ & 3.623e-08 & 1.929e-07 & 5.686e-07 & 1.277e-06 \\
$p=15$ &  8.778e-12 & 9.668e-11 & 4.670e-10 & 1.575e-09      
\end{tabular}
}
\qquad
  \subfloat[Optimized metrics.]{
\begin{tabular}{l|c|c|c|c}
  & \multicolumn{4}{c}{$\mathcal{T}$}  \\ 
\cline{1-5}
$\eta$          & $0.25$ & $0.5$ & $0.75$  & $1.0$  \\ 
\hline
$p=1 $ & 3.038e-03 & 3.038e-03 & 3.038e-03 & 3.038e-03 \\
$p=2 $ & 4.350e-04 & 5.648e-04 & 7.632e-04 & 1.022e-03 \\
$p=3 $ & 7.775e-05 & 1.489e-04 & 2.584e-04 & 4.058e-04 \\
$p=4 $ & 2.859e-05 & 6.747e-05 & 1.265e-04 & 2.094e-04 \\
$p=7 $ & 4.951e-07 & 2.121e-06 & 5.425e-06 & 1.102e-05 \\
$p=9 $ & 3.506e-08 & 1.910e-07 & 5.669e-07 & 1.278e-06 \\
$p=15$ &  8.735e-12 & 9.668e-11 & 4.670e-10 & 1.575e-09           
\end{tabular}
}
  \caption{$L^2$ error norm of the temperature: 
  propagation of the isentropic vortex in a cubic domain with geometrically high order perturbed cells.}
\label{tbl:ivortex_cube_t}
\end{table}

\pagebreak
\subsection{Cubic domain with a spherical hole and geometrically high order cells}

\begin{table}[htbp!]
\begin{center}
\begin{tabular}{l|c|c|c|c|c}
          & $\rho$ & $\mathcal{U}_{1}$ & $\mathcal{U}_{2}$  & $\mathcal{U}_{3}$ & $\mathcal{T}$  \\ 
\hline
 $ p=2$ & 1.264e-03  & 6.661e-03  & 6.661e-03  & 3.507e-03   & 5.725e-04  \\
 $ p=3$ & 1.293e-04  & 6.434e-04  & 6.434e-04  & 6.700e-04   & 5.397e-05  \\
 $ p=4$ & 2.193e-05  & 1.167e-04  & 1.167e-04  & 6.360e-05   & 1.208e-05  \\
 $ p=7$ & 5.360e-08  & 5.091e-07  & 5.091e-07  & 2.593e-07   & 2.957e-08  \\
 $ p=9$ & 1.608e-09  & 1.762e-08  & 1.762e-08  & 7.390e-09   & 7.782e-10         
\end{tabular}
\caption{$L^2$-norm of the errors for the Thomas and Lombard metrics \cite{thomas_gcl_1979}: 
  propagation of the isentropic vortex in a cubic domain with a spherical hole and geometrically high order cells.}
\label{tbl:ivortex_cube_err}
\end{center}
\end{table}

\begin{table}[htbp!]
\begin{center}
\begin{tabular}{l|c|c|c|c|c}
          & $\rho$ & $\mathcal{U}_{1}$ & $\mathcal{U}_{2}$  & $\mathcal{U}_{3}$ & $\mathcal{T}$  \\ 
\hline
 $ p=2$ & 1.155e-03  & 6.665e-03  & 6.665e-03   & 3.079e-03  & 5.236e-04 \\
 $ p=3$ & 1.270e-04  & 6.311e-04  & 6.311e-04   & 6.653e-04  & 5.268e-05 \\
 $ p=4$ & 2.074e-05  & 1.164e-04  & 1.164e-04   & 6.187e-05  & 1.189e-05 \\
 $ p=7$ & 5.149e-08  & 5.019e-07  & 5.019e-07   & 2.496e-07  & 2.889e-08 \\
 $ p=9$ & 1.562e-09  & 1.742e-08  & 1.742e-08   & 7.268e-09  & 7.614e-10      
\end{tabular}
\caption{$L^2$-norm of the errors for the optimized metrics: 
  propagation of the isentropic vortex in a cubic domain with a spherical hole and geometrically high order cells.}
\label{tbl:ivortex_cube_sphere_err}
\end{center}
\end{table}

\pagebreak
\section{Numerical results: Viscous shock}
\subsection{Cubic domain with geometrically high order perturbed cells}

\begin{table}
\centering
  \subfloat[Thomas and Lombard metrics \cite{thomas_gcl_1979}.]{
\begin{tabular}{l|c|c|c|c}
 & \multicolumn{4}{c}{$\rho$}  \\ 
\cline{1-5}
$\eta$          & $0.25$ & $0.5$ & $0.75$  & $1.0$  \\ 
\hline
$p=1 $ & 3.795e-02 & 3.795e-02 & 3.795e-02 & 3.795e-02 \\
$p=2 $ & 1.673e-02 & 2.486e-02 & 3.384e-02 & 4.267e-02 \\
$p=3 $ & 4.121e-03 & 6.915e-03 & 9.867e-03 & 1.282e-02 \\
$p=4 $ & 1.956e-03 & 3.350e-03 & 5.118e-03 & 7.082e-03 \\
$p=7 $ & 1.244e-04 & 2.699e-04 & 4.902e-04 & 7.852e-04 \\
$p=9 $ & 2.216e-05 & 5.799e-05 & 1.192e-04 & 2.118e-04 \\
$p=15$ &  1.393e-07 & 6.844e-07 & 2.249e-06 & 5.832e-06      
\end{tabular}
}
\qquad
  \subfloat[Optimized metrics.]{
\begin{tabular}{l|c|c|c|c}
 & \multicolumn{4}{c}{$\rho$}  \\ 
\cline{1-5}
$\eta$          & $0.25$ & $0.5$ & $0.75$  & $1.0$  \\ 
\hline
$p=1 $ & 3.795e-02 & 3.795e-02 & 3.795e-02 & 3.795e-02 \\
$p=2 $ & 1.380e-02 & 1.618e-02 & 1.906e-02 & 2.187e-02 \\
$p=3 $ & 4.014e-03 & 6.554e-03 & 9.156e-03 & 1.168e-02 \\
$p=4 $ & 1.863e-03 & 3.125e-03 & 4.840e-03 & 6.862e-03 \\
$p=7 $ & 1.249e-04 & 2.713e-04 & 4.919e-04 & 7.860e-04 \\
$p=9 $ & 2.217e-05 & 5.800e-05 & 1.192e-04 & 2.117e-04  \\
$p=15$ &  1.393e-07 & 6.844e-07 & 2.249e-06 & 5.832e-06            
\end{tabular}
}
\caption{$L^2$ error norm of the density: 
  propagation of the viscous shock in a cubic domain with geometrically high order perturbed cells.}
\label{tbl:vs_cube_rho_err}
\end{table}

\begin{table}
\centering
  \subfloat[Thomas and Lombard metrics \cite{thomas_gcl_1979}.]{
\begin{tabular}{l|c|c|c|c}
  & \multicolumn{4}{c}{$\Um{1}$}  \\ 
\cline{1-5}
$\eta$          & $0.25$ & $0.5$ & $0.75$  & $1.0$  \\ 
\hline
$p=1 $ & 1.030e-02 & 1.030e-02 & 1.030e-02 & 1.030e-02 \\
$p=2 $ & 2.406e-03 & 3.587e-03 & 4.915e-03 & 6.335e-03 \\
$p=3 $ & 6.446e-04 & 9.659e-04 & 1.346e-03 & 1.769e-03 \\
$p=4 $ & 2.330e-04 & 4.228e-04 & 6.496e-04 & 9.159e-04 \\
$p=7 $ & 1.023e-05 & 2.333e-05 & 4.458e-05 & 7.568e-05 \\
$p=9 $ & 1.529e-06 & 4.296e-06 & 9.478e-06 & 1.777e-05  \\
$p=15$ &  6.950e-09 & 3.633e-08 & 1.228e-07 & 3.236e-07       
\end{tabular}
}
\qquad
  \subfloat[Optimized metrics.]{
\begin{tabular}{l|c|c|c|c}
  & \multicolumn{4}{c}{$\Um{1}$}  \\ 
\cline{1-5}
$\eta$          & $0.25$ & $0.5$ & $0.75$  & $1.0$  \\ 
\hline
$p=1 $ & 1.030e-02 & 1.030e-02 & 1.030e-02 & 1.030e-02 \\
$p=2 $ & 2.027e-03 & 2.530e-03 & 3.180e-03 & 3.904e-03 \\
$p=3 $ & 6.236e-04 & 9.166e-04 & 1.271e-03 & 1.658e-03 \\
$p=4 $ & 2.122e-04 & 3.773e-04 & 5.832e-04 & 8.344e-04 \\
$p=7 $ & 1.031e-05 & 2.366e-05 & 4.530e-05 & 7.687e-05 \\
$p=9 $ & 1.531e-06 & 4.306e-06 & 9.496e-06 & 1.779e-05  \\
$p=15$ &  6.950e-09 & 3.633e-08 & 1.228e-07 & 3.236e-07             
\end{tabular}
}
  \caption{$L^2$ error norm of the velocity component in $\xm{1}$: 
  propagation of the viscous shock in a cubic domain with geometrically high order perturbed cells.}
\label{tbl:vs_cube_u1_err}
\end{table}

\begin{table}
\centering
  \subfloat[Thomas and Lombard metrics \cite{thomas_gcl_1979}.]{
\begin{tabular}{l|c|c|c|c}
  & \multicolumn{4}{c}{$\Um{2}$}  \\ 
\cline{1-5}
$\eta$          & $0.25$ & $0.5$ & $0.75$  & $1.0$  \\ 
\hline
$p=1 $ & 1.030e-02 & 1.030e-02 & 1.030e-02 & 1.030e-02 \\
$p=2 $ & 2.406e-03 & 3.587e-03 & 4.915e-03 & 6.335e-03 \\
$p=3 $ & 6.446e-04 & 9.659e-04 & 1.346e-03 & 1.769e-03 \\
$p=4 $ & 2.330e-04 & 4.228e-04 & 6.496e-04 & 9.159e-04 \\
$p=7 $ & 1.023e-05 & 2.333e-05 & 4.458e-05 & 7.568e-05 \\
$p=9 $ & 1.529e-06 & 4.296e-06 & 9.478e-06 & 1.777e-05  \\
$p=15$ &  6.950e-09 & 3.633e-08 & 1.228e-07 & 3.236e-07          
\end{tabular}
}
\qquad
  \subfloat[Optimized metrics.]{
\begin{tabular}{l|c|c|c|c}
  & \multicolumn{4}{c}{$\Um{2}$}  \\ 
\cline{1-5}
$\eta$          & $0.25$ & $0.5$ & $0.75$  & $1.0$  \\ 
\hline
$p=1 $ & 1.030e-02 & 1.030e-02 & 1.030e-02 & 1.030e-02 \\
$p=2 $ & 2.027e-03 & 2.530e-03 & 3.180e-03 & 3.904e-03 \\
$p=3 $ & 6.236e-04 & 9.166e-04 & 1.271e-03 & 1.658e-03 \\
$p=4 $ & 2.122e-04 & 3.773e-04 & 5.832e-04 & 8.344e-04 \\
$p=7 $ & 1.031e-05 & 2.366e-05 & 4.530e-05 & 7.687e-05 \\
$p=9 $ & 1.531e-06 & 4.306e-06 & 9.496e-06 & 1.779e-05 \\
$p=15$ &  6.950e-09 & 3.633e-08 & 1.228e-07 & 3.236e-07          
\end{tabular}
}
  \caption{$L^2$ error norm of the velocity component in $\xm{2}$: 
  propagation of the viscous shock in a cubic domain with geometrically high order perturbed cells.}
\label{tbl:vs_cube_u2_err}
\end{table}

\begin{table}
\centering
  \subfloat[Thomas and Lombard metrics \cite{thomas_gcl_1979}.]{
\begin{tabular}{l|c|c|c|c}
  & \multicolumn{4}{c}{$\Um{3}$}  \\ 
\cline{1-5}
$\eta$          & $0.25$ & $0.5$ & $0.75$  & $1.0$  \\ 
\hline
$p=1 $ & 1.030e-02 & 1.030e-02 & 1.030e-02 & 1.030e-02 \\
$p=2 $ & 2.406e-03 & 3.587e-03 & 4.915e-03 & 6.335e-03 \\
$p=3 $ & 6.446e-04 & 9.659e-04 & 1.346e-03 & 1.769e-03 \\
$p=4 $ & 2.330e-04 & 4.228e-04 & 6.496e-04 & 9.159e-04 \\
$p=7 $ & 1.023e-05 & 2.333e-05 & 4.458e-05 & 7.568e-05\\ 
$p=9 $ & 1.529e-06 & 4.296e-06 & 9.478e-06 & 1.777e-05 \\
$p=15$ &  6.950e-09 & 3.633e-08 & 1.228e-07 & 3.236e-07       
\end{tabular}
}
\qquad
  \subfloat[Optimized metrics.]{
\begin{tabular}{l|c|c|c|c}
  & \multicolumn{4}{c}{$\Um{3}$}  \\ 
\cline{1-5}
$\eta$          & $0.25$ & $0.5$ & $0.75$  & $1.0$  \\ 
\hline
$p=1 $ & 1.030e-02 & 1.030e-02 & 1.030e-02 & 1.030e-02 \\
$p=2 $ & 2.027e-03 & 2.530e-03 & 3.180e-03 & 3.904e-03 \\
$p=3 $ & 6.236e-04 & 9.166e-04 & 1.271e-03 & 1.658e-03 \\
$p=4 $ & 2.122e-04 & 3.773e-04 & 5.832e-04 & 8.344e-04 \\
$p=7 $ & 1.031e-05 & 2.366e-05 & 4.530e-05 & 7.687e-05 \\
$p=9 $ & 1.531e-06 & 4.306e-06 & 9.496e-06 & 1.779e-05 \\
$p=15$ &  6.950e-09 & 3.633e-08 & 1.228e-07 & 3.236e-07           
\end{tabular}
}
  \caption{$L^2$ error norm of the velocity component in $\xm{3}$: 
  propagation of the viscous shock in a cubic domain with geometrically high order perturbed cells.}
\label{tbl:vs_cube_u3_err}
\end{table}

\begin{table}
\centering
  \subfloat[Thomas and Lombard metrics \cite{thomas_gcl_1979}.]{
\begin{tabular}{l|c|c|c|c}
  & \multicolumn{4}{c}{$\mathcal{T}$}  \\ 
\cline{1-5}
$\eta$          & $0.25$ & $0.5$ & $0.75$  & $1.0$  \\ 
\hline
$p=1 $ & 2.864e-02 & 2.864e-02 & 2.864e-02 & 2.864e-02 \\
$p=2 $ & 5.072e-03 & 8.012e-03 & 1.142e-02 & 1.530e-02 \\
$p=3 $ & 1.190e-03 & 1.982e-03 & 2.942e-03 & 4.102e-03 \\
$p=4 $ & 4.457e-04 & 8.996e-04 & 1.441e-03 & 2.090e-03 \\
$p=7 $ & 1.284e-05 & 3.472e-05 & 7.484e-05 & 1.366e-04\\ 
$p=9 $ & 1.576e-06 & 5.798e-06 & 1.437e-05 & 2.818e-05 \\
$p=15$ &  5.170e-09 & 3.352e-08 & 1.342e-07 & 4.256e-07      
\end{tabular}
}
\qquad
  \subfloat[Optimized metrics.]{
\begin{tabular}{l|c|c|c|c}
  & \multicolumn{4}{c}{$\mathcal{T}$}  \\ 
\cline{1-5}
$\eta$          & $0.25$ & $0.5$ & $0.75$  & $1.0$  \\ 
\hline
$p=1 $ & 2.864e-02 & 2.864e-02 & 2.864e-02 & 2.864e-02 \\
$p=2 $ & 4.133e-03 & 5.462e-03 & 7.282e-03 & 9.373e-03 \\
$p=3 $ & 1.077e-03 & 1.677e-03 & 2.400e-03 & 3.250e-03 \\
$p=4 $ & 3.704e-04 & 7.448e-04 & 1.218e-03 & 1.801e-03 \\
$p=7 $ & 1.289e-05 & 3.435e-05 & 7.336e-05 & 1.336e-04\\ 
$p=9 $ & 1.573e-06 & 5.764e-06 & 1.428e-05 & 2.803e-05 \\
$p=15$ &  5.170e-09 & 3.352e-08 & 1.342e-07 & 4.256e-07          
\end{tabular}
}
  \caption{$L^2$ error norm of the temperature: 
  propagation of the viscous shock in a cubic domain with geometrically high order perturbed cells.}
\label{tbl:vs_cube_t_err}
\end{table}

\pagebreak
\subsection{Cubic domain with a spherical hole and geometrically high order cells}

\begin{table}[htbp!]
\begin{center}
\begin{tabular}{l|c|c|c|c|c}
          & $\rho$ & $\mathcal{U}_{1}$ & $\mathcal{U}_{2}$  & $\mathcal{U}_{3}$ & $\mathcal{T}$  \\ 
\hline
$p=2 $ & 1.184e-02      & 1.519e-03     & 1.519e-03     & 1.519e-03           & 4.269e-03\\
$p=3 $ & 4.587e-03      & 4.340e-04     & 4.340e-04     & 4.340e-04           & 1.057e-03\\
$p=4 $ & 1.820e-03      & 1.340e-04     & 1.340e-04     & 1.340e-04           & 3.222e-04\\
$p=7 $ & 1.370e-04      & 9.784e-06     & 9.784e-06     & 9.784e-06           & 1.054e-05\\
$p=9 $ & 2.176e-05      & 1.854e-06     & 1.854e-06     & 1.854e-06           & 1.729e-06         
\end{tabular}
\caption{$L^2$-norm of the errors for the Thomas and Lombard metrics \cite{thomas_gcl_1979}: 
  propagation of the viscous shock in a cubic domain with a spherical hole and geometrically high order cells.}
\label{tbl:vs_cube_sphere_err_tl}
\end{center}
\end{table}

\begin{table}[htbp!]
\begin{center}
\begin{tabular}{l|c|c|c|c|c}
          & $\rho$ & $\mathcal{U}_{1}$ & $\mathcal{U}_{2}$  & $\mathcal{U}_{3}$ & $\mathcal{T}$  \\ 
\hline
$p=2 $ & 1.169e-02      & 1.548e-03     & 1.548e-03     & 1.548e-03           & 4.076e-03 \\
$p=3 $ & 4.586e-03      & 4.352e-04     & 4.352e-04     & 4.352e-04           & 1.062e-03 \\
$p=4 $ & 1.817e-03      & 1.340e-04     & 1.340e-04     & 1.340e-04           & 3.205e-04 \\
$p=7 $ & 1.370e-04      & 9.784e-06     & 9.784e-06     & 9.784e-06           & 1.054e-05 \\
$p=9 $ & 2.176e-05      & 1.854e-06     & 1.854e-06     & 1.854e-06           & 1.729e-06       
\end{tabular}
\caption{$L^2$-norm of the errors for the optimized metrics: 
  propagation of the viscous shock in a cubic domain with a spherical hole and geometrically high order cells.}
\label{tbl:vs_cube_sphere_err_opt}
\end{center}
\end{table}

\pagebreak

\end{document}